\documentclass[preprint,11pt]{article}
\usepackage{amssymb}
\usepackage{mathrsfs}
\usepackage{amsfonts} 
\usepackage{indentfirst}
\usepackage{colortbl,dcolumn}
\usepackage{color}
\usepackage{amsmath}
\usepackage{psfrag}
\usepackage{booktabs}
\usepackage{cite}
\usepackage[hidelinks]{hyperref}
\usepackage[affil-it]{authblk}
\usepackage[english]{babel}
\usepackage{blindtext}
\makeatletter
\newcommand\tabcaption{\def\@captype{table}\caption}
\newcommand\figcaption{\def\@captype{figure}\caption}
\makeatother
\usepackage{graphicx} 
\graphicspath{{figs/}}
\usepackage{amsthm}
\numberwithin{equation}{section}
\usepackage[top=1in, bottom=1in, left=0.8in, right=0.8in]{geometry}
\newcommand{\R}{\mathbb{R}}
\newcommand{\N}{\mathbb{N}}
\newcommand{\E}{\mathbb{E}}
\renewcommand{\P}{\mathbb{P}}

\newcommand{\dd}{\text{d}}

\newtheorem{thm}{Theorem}[section]
\newtheorem{defn}[thm]{Definition}
\newtheorem{lem}[thm]{Lemma}
\newtheorem{exm}[thm]{Example}
\newtheorem{prop}[thm]{Proposition}
\newtheorem{cor}[thm]{Corollary}

\newtheorem{example}[thm]{Example}
\newtheorem{assumption}[thm]{Assumption}

\allowdisplaybreaks
\begin{document}
\title{Weak error analysis for strong approximation schemes of SDEs with super-linear coefficients
\footnotemark[2] \footnotetext[2]{This work was supported by Natural Science Foundation of China (12071488, 11971488),
                Natural Science Foundation of Hunan Province (2020JJ2040) 
                and the innovative project of graduate students of  Central South University  (2021zzts0040). \\
                E-mail addresses: x.j.wang7@csu.edu.cn, yzhao15@wpi.edu,
                zzhang7@wpi.edu
                }}
\author[a]{Xiaojie Wang}
\author[a,b]{Yuying Zhao \thanks{corresponding author}}
\author[b]{Zhongqiang Zhang}
\affil[a]{
    School of Mathematics and Statistics, HNP-LAMA, Central South University, Changsha, Hunan, P. R. China}
\affil[b]{
    Department of Mathematical Sciences, Worcester Polytechnic Institute, Worcester, MA 01609 USA}
\maketitle	
\begin{abstract}
We present an error analysis of weak convergence of one-step numerical schemes for stochastic differential equations (SDEs) with super-linearly growing coefficients. Following  Milstein's weak error analysis on the one-step approximation of SDEs, we prove  a general conclusion on weak convergence of the one-step discretization of the SDEs mentioned above. As applications, we show the  weak convergence rates for several numerical schemes of half-order strong convergence, such as tamed and balanced schemes. Numerical examples are presented to  verify our theoretical analysis.
	
\end{abstract}
\tableofcontents

\section{Introduction}
\noindent
In the past three decades, 
numerical methods for stochastic differential equations (SDEs) from various applications have been extensively explored.
When the coefficients of SDEs are globally Lipschitz continuous, numerical methods and their analysis 
are well documented in  
\cite{Milstein2004stochastic,Kloeden1992stochastic,Platen10numerical},
see particularly \cite{Milstein1978method, Abdulle2012high, Abdulle2013weak, Rossler2006rooted, Debrabant2009b, Komori2017weak,Buckwar2005weak,talay1990expansion,Lyons2004cubature,Kohatsu2001weak,Jourdain2011review,Pageesintroduction2018} and references therein
for further analysis of weak approximation schemes.
When the coefficients of SDEs are non-globally Lipschitz  continuous with super-linear growth, 
it is shown that the commonly used Euler-Maruyama  method fails to converge in the strong and weak sense (see  \cite{Milstein2005Numerical}, \cite{Higham2007almost},
\cite{mattingly2002ergodicity}
and
\cite{Hutzenthaler2011Strong}).
 Since then many papers are focusing on SDEs with super-linearly growing coefficients, such as 
\cite{Higham2007almost,Hutzenthaler2011Strong, Tretyakov2013fundamental,Hutzenthaler2012Strong,Hutzenthaler2015Numerical,Sabanis2016euler,
Mao2015truncated,Fang2020adaptive,Kelly2017adaptive,Beyn2016stochastic,Brehier2020weak,brehier2020approximation,sabanis2013note} and references therein. 

Most works focus on strong approximations of SDEs with superlinearly growing coefficients. 
There are fewer works on the weak approximation of such SDEs.  
In \cite{Milstein2005Numerical}, the authors
introduced a new concept of 
discarding the approximate trajectories that leave a sufficiently large sphere  $ S_R:=\{x:|x|<R\} $. 
This strategy allows {\color{black}the authors} to use any usual weak scheme for solving a {\color{black}broad} class of SDEs with non-globally Lipschitz coefficients,
resulting in arbitrarily small weak errors by increasing the radius of the sphere. However, the weak convergence rates are not explicitly revealed.
{\color{black} In \cite{Mora2017stable}, the authors devised
a new kind of direction and norm decomposition methods for SDEs with superlinearly growing coefficients and showed 
a linear rate of weak convergence and one-half strong order of convergence under mild assumptions.}
The authors of \cite{bossy2021weak} proposed an exponential Euler scheme for one-dimensional SDE 
with superlinearly growing coefficients and recovered a rate of weak convergence of order one for the scheme.
{\color{black} In \cite{brehier2020approximation}, the author analyzed the long-time error estimates of an explicit tamed Euler scheme applied to a class of SDEs driven by additive noise, under a one-sided Lipschitz continuity and superlinear growth condition.

In addition, the topic of weak error analysis for SDEs with non-smooth coefficients or non-smooth test function has and remains an active field of research which have made tremendous progress over the past years.
For example, we observe that there are several works on
weak convergence of  non-globally Lipschitz coefficients of \textit{sublinear} growth  (e.g, the Heston model), e.g., \cite{mickel2021weak,zheng2017weak,altmayer2017discretising}.
There are also many works on weak approximations of SDEs with 
singular coefficients or test functions,
see, e.g., \cite{Bally1996law,bao2020convergence,suo2022weak,shao2018weak,kohatsu2017weak,ngo2018approximation,Konakov2002edgeworth,Konakov2017weak,Guyon2006euler,Bencheikh2022convergence,Jourdain2021convergence,Yan2002euler}, where the classical scheme such as the Euler-Maruyama scheme still provides an efficient and quantifiable numerical approximation. 
The common ground of these references is roughly to establish properties of the transition semi-groups related to the considered SDE. 
}

{\color{black}In this work, we  develop a weak error analysis for general one-step approximation schemes of SDEs with non-globally Lipschitz coefficients of superlinear growth (see Assumption \ref{ass:coefficient_function_assumption}). 
Specifically, we first {\color{black}revisit} the general weak convergence theorem for general one-step numerical approximations and reformulate it under mild conditions (see Theorem \ref{thm:fundamental_weak_convergence}) following Milstein's framework, see e.g., Theorem 2 of \cite{Milstein1985Weak} and Theorem 
2.1 in Chapter 2 of \cite{Milstein2004stochastic}.
We then establish some explicit conditions on the coefficients of SDEs to guarantee the aforementioned assumptions in the general convergence theorem.
These assumptions are discussed in 
Section \ref{section:setting_gen_the} 
and we then derive a corollary on weak convergence of numerical methods for SDEs with one-sided Lipschitz continuous and superlinearly growing coefficients  (Corollary \ref{cor-fundamental-weak-convergence}).}
As applications of the corollary, we prove a first-order (or less) weak convergence for various modified Euler approximations (Theorem \ref{thm:MEM-global-conver-rate}), which are based on strong convergence,  
and the backward Euler method (Theorem \ref{thm:BEM-global-conver-rate}). 
We summarize our main contributions:
\begin{itemize}
	\item
{\color{black} The weak convergence theorem is revised for general one-step numerical approximations
under mild conditions (see Theorem \ref{thm:fundamental_weak_convergence}).
Based on the theorem, a more specific conclusion on weak convergence of numerical methods is obtained for SDEs with one-sided Lipschitz continuous and superlinearly growing coefficients
(see Corollary \ref{cor-fundamental-weak-convergence})
 }
	\item 
Weak convergence rates for a class of  modified Euler methods (Theorem \ref{thm:MEM-global-conver-rate}) and the backward Euler method
(Theorem \ref{thm:BEM-global-conver-rate}) are established, which fills some notable gaps in literature on the weak error analysis of strong approximating schemes.
\end{itemize}	 
\textbf{How to apply our results.} 
Roughly speaking, when the conditions in Assumption \ref{ass:coefficient_function_assumption} are satisfied, Corollary \ref{cor-fundamental-weak-convergence} asserts that if moments of the numerical solution are bounded and the scheme was proved to be convergent with order $p+1$ in the one-step error (Assumption \ref{ass:local-error}), then the scheme has the weak convergence order $p$. 
For example, the boundedness of moments for the modified Euler schemes is proved under Assumption \ref{con:modified_EM_dri_diff} on the bounds of modified drift and  the diffusion coefficients (Theorem \ref{thm:bound_moment_modi_EM}). Then, with Assumption \ref{ass:error-drift-diff-modi-EM} on the errors of modifying drift and  the diffusion coefficients 
we can obtain from Theorem 
\ref{thm:MEM-global-conver-rate} that
the modified Euler schemes \eqref{eq:modified_EM_scheme} converge with weak order $ p = 1$ or smaller, depending on the choices of modifications of the coefficients; see examples in Section \ref{ssec:example-mod-euler}.

%

The rest of this article is structured as follows.
{\color{black}
In the next section,  we revisit the general convergence theorem on weak convergence of one-step schemes for SDEs with 
non-globally Lipschitz continuous coefficients.}
In Section \ref{section:corollaries-weak-convergence}, we present explicit conditions on coefficients of SDEs, where the assumptions are satisfied in the general convergence theorem. 
To prove a first-order weak convergence, we compare the considered schemes with  the Euler-Maruyama method and its one-step error 
 in Section \ref{sec:one-step-euler}. The one-step error of  the Euler-Maruyama method is proved in Appendix \ref{appendix:proof_EM_one_step_error}.
In Sections \ref{sec:weak-rate-MEM} and \ref{sec:weak-rate-BEM}, 
we apply the corollary on weak convergence obtained in Section \ref{section:corollaries-weak-convergence} to derive weak convergence rates for some explicit and implicit  numerical methods.
We present some numerical results in Section \ref{sec:numer-results}. 
\section{The general weak convergence theorem revisited} 
\label{section:setting_gen_the}
\noindent
In this section, we introduce notations and basic assumptions on stochastic differential equations. We also introduce the one-step approximation and
revisit a weak convergence theorem for it.
\subsection{Notations}
Throughout this paper, the following notation
is frequently used.
Let $ | \cdot | $ and  
$ \left\langle \cdot , \cdot \right\rangle $
be the Euclidean norm and the inner product of vectors
in $ \R^d $, respectively.
By $ A^T $ we denote the transpose of vector or
matrix $ A $. 
Given a matrix $ A $,  we use
$ \| A \| : = \sqrt{ \text {trace} (A^TA) } $
to denote the trace norm of $A$.
On a probability space $ ( \Omega, \mathcal{ F }, \P ) $, we use $ \E $ to mean expectation and
$ L^r(\Omega; \R^{d \times m}) $, $ r \in \mathbb{N} $,
to denote the family of 
$ \R^{d \times m} $-valued variables with
the norm defined by 
$ \| \xi \|_{L^r(\Omega;\R^{d \times m})} 
=(\E [|\xi|^r])^{\frac{1}{r}} < \infty $.
Next, we give some notation for partial derivatives of functions, frequently used throughout this paper.
	A vector $\alpha=\left(\alpha_{1}, \ldots, \alpha_{d}\right)$ 
	is called a multiindex of order
	$
	|\alpha|=\alpha_{1}+\cdots+\alpha_{d},
	$
	where each component $\alpha_{i} \geq 0$  is a nonnegative integer.
	For a multiindex $\alpha$, we define the partial derivatives of $v \colon \R^d \rightarrow \R^l $ as
        \begin{equation}
	D^{\alpha} v(x)
	:=
	\tfrac{\partial^{|\alpha|} 
		v(x)}{\partial x_{1}^{\alpha_{1}} 
		\cdots \partial x_{d}^{\alpha_{d}}}
	=\partial_{x_{1}}^{\alpha_{1}} 
	\cdots \partial_{x_{d}}^{\alpha_{d}} v (x) .
	\end{equation}
	For a nonnegative integer $k$, we use
	$
	D^{k} v ( x ) := \left\{D^{\alpha} v ( x ) 
	: | \alpha | =k \right\}
	$  
	to denote the set of all partial derivatives of order $k$.  Moreover, we define  
	\begin{equation}
	\left|D^{k}  v ( x ) \right|
	=\Big( \sum_{|\alpha|=k}\left|D^{\alpha} 
	v (x)\right|^{2} \Big)^{1 / 2}.
	\end{equation}
%
%
%
%
   We also use 
   $  D v := 
             (\tfrac{\partial v}{ \partial x_{1} },
                   \tfrac{\partial v}{ \partial x_{2} },
                    \cdots
                   \tfrac{\partial v}{ \partial x_{d} })
                   $  (gradient  vector)  and $
        D^2  v :=
                 \big( \frac{\partial ^{2} v}
                 {\partial x_{i} \partial x_{j}}\big)_{d \times d}
             $ (Hessian matrix). 

Next, we introduce a class of functions with superlinear growth.
\begin{defn}[polynomial growth]
 \label{defi:function_class_H}
	A function $ \psi \colon \mathbb{R}^{d} \rightarrow \mathbb{R} $
	is said to belong to the class $ \mathbb{H} $, written as
	$ \psi \in \mathbb{H} $,
	if there exist constants $ L > 0, l > 0 $ such that the following inequality holds 
	\begin{equation}
		\label{eq:function-class-H}
		|\psi( x )| \leq L ( 1 + | x |^{ l } ).
	\end{equation}
For $k \in \mathbb{N}$ 
	we also use $\mathbb{H}^{k}$ to denote a subset of $\mathbb{H}$, consisting of $k$-times continuously 
	differentiable functions 
	which, together with its partial derivatives 
	up to and including order $k$, 
	belong to $\mathbb{H}$. 
\end{defn}
{\color{black}The above definition is quoted from \cite[Definition 1.1 in Chapter 2]{Milstein2004stochastic}}.
For random functions, we introduce a concept of mean-square 
differentiability as follows.
\begin{defn}[mean-square differentiable]
	\label{defi:mean_square_differential}
Let 
   $ \Psi \colon 
     \Omega \times \mathbb{R}^{d} \rightarrow
       \mathbb{R} $
   and 
   	$  \phi_{i} \colon \Omega \times \mathbb{R}^{d} 
   	 \rightarrow \mathbb{R} $
be random functions satisfying
   \begin{equation}
     \lim_{\tau \rightarrow 0 }
        \E 
        \Big [ \Big|
         \frac{1}{\tau}
         \big [
          \Psi( x + \tau e_{i})
        -
          \Psi( x ) \big]
         -
          \phi_i(x) 
          \Big|^2
         \Big]
        = 0,
        \quad
        \forall
i \in \{ 1,2,\cdots, d \},
   \end{equation}
where $e_i \in \R^d $ is a unit vector in $\R^d$, with the $i$-th element being $1$. 
Then $ \Psi  $ is called to be mean-square differentiable, with
$ \phi = (\phi_{1}, \phi_{2}, \cdots, \phi_{d} ) $
being the derivative (in the mean-square differentiable sense)
of $ \Psi  $ and  we also write  $ \partial_{(i)}\Psi = \phi_{i} $.
\end{defn}
The above two definitions can be generalized to vector-valued functions in a component-wise manner.
%
\subsection{Stochastic differential equations}
Let $ ( \Omega, \mathcal{ F }, \P ) $ 
be a complete probability space with a 
filtration 
$ \{ \mathcal{ F }_t \}_{ t \in [0,T] } $  ($ T > 0 $),
satisfying the usual conditions.
Let
$ d, m \in \N $, and 
$ W: [0,T] \times \Omega 
\rightarrow \R^m $
stands for 
standard Brownian motions
with respect to 
$ \{ \mathcal{ F }_t \}_{ t \in [0,T] } $.
We consider the following autonomous SDEs in the It\^o form of
\begin{equation}
	\label{eq:Problem_SDE}
	\left\{
	\begin{aligned}
		\dd X(t) & = f ( X(t) ) \, \dd t 
		+ g ( X(t) ) \, \dd W(t),
		 \quad t \in (0,T],
		  \\
		X(0) & = X_{0} \in \mathbb{R}^d,
	\end{aligned}\right.
\end{equation}
where $f = (f^{1},f^{2},\cdots,f^{d})^T \colon \mathbb{R}^d \rightarrow \mathbb{R}^d$
is the drift coefficient and $f^{i} \colon \R^d \rightarrow \R$, and
$ g= (g^{i,j})_{d \times m} \colon \mathbb{R}^d \rightarrow \mathbb{R}^{d\times m} $
is the diffusion coefficient. The function $g$ is also written as 
$ g  
= (g^{1},g^{2},\cdots,g^{m})
=
 (g_{1}^T,g_{2}^T,\cdots,g_{d}^T)^T 
$, where 
$g^j  \colon \R^d \rightarrow \R^d $  is a $d \times 1$ vector function
and
$g_i  \colon \R^d \rightarrow \R^m $  is a $1 \times {m} $ vector function.
The initial data 
$ X_0 \in \R^d $ is assumed to be deterministic for simplicity. 
\begin{assumption}[Well-posedness of a strong solution]
	\label{ass:SDEs-wellposed}
	{\color{black}For any $X_0 \in  \R^d$, the SDE \eqref{eq:Problem_SDE} admits 
	a unique pathwise solution}  
	$ X:[0,T] \times \Omega \rightarrow \R^d $, given by
	\begin{equation}\label{eq:solution-of-SDE}
	X(t) =
	X_0 + \int_{0}^{t} f(X(s)) \, \dd s
	+ \int_{0}^{t} g( X(s) )  \, \dd W(s),
	\quad
	a.s..
	\end{equation} 
\end{assumption}
Let $X(t,x;s)$ denotes the solution to \eqref{eq:Problem_SDE} at $s$, starting from the initial value $ x$ at $t$, namely,
\begin{equation}
\label{eq:SDE-X(t,x;T)}
X(t,x;s) = x 
+ 
\int_t^s f( X ( t, x; r ) ) \, \dd r 
+
\int_t^s g( X ( t, x; r ) ) \, \dd W ( r ),
\quad
0 \leq t \leq s \leq T.
\end{equation}
\subsection{The general weak convergence theorem}
To approximate solutions to  
\eqref{eq:Problem_SDE}, we construct a uniform mesh on
$ [0,T] $ with  $ h = \frac{T}{N} $ being the 
step size,
for any $ N \in \mathbb{N} $.
Here we use $X_{t,x}(s)$ or $X(t,x;s)$ for $0 \leq s \leq T$ to denote the solution to \eqref{eq:Problem_SDE} at the time instant $s$ with the initial condition
$X_{t,x}(t)=X(t,x;t)=x$.
Also, $ X(t) := X(0,X_0;t) $.
For $ x \in \R^d $, $ t \in [0,T] $, $ h > 0 $,
$ 0 < t+h \leq T $,
we introduce the one-step approximation
$ Y(t,x;t+h) $ 
in the form of
\begin{equation}
	\label{sq:numerical-solution-sequence}
	Y(t,x;t+h)= x + \Phi( t,x,h;\xi_t ),
\end{equation}
where $ \xi_t $ is a random variable  (in general, a vector) defined on $ ( \Omega, \mathcal{ F }, \P ) $,
having moments of a sufficiently high order,
and $ \Phi $ is a function from
$ [0,T] \times \R^d \times (0,T] \times \R^m $
to $ \R^d $.
Similarly, by $Y(t,x;t+h)$ or 
$Y_{t,x}(t+h)$ we denote an approximation of
the solution at $t+h$ 
with initial value 
$ Y(t,x;t) = Y_{t,x}(t) = x $.
Then we use 
$ Y_{n+1} = Y(t_n,Y_n;t_{n+1}) $
to recurrently construct numerical approximations
$ \{Y_n\}_{0 \leq n \leq N} $
on the uniform mesh grid 
$ \{t_n = nh, n = 0, 1, \cdots, N\} $,
given by
\begin{equation}
	\label{eq:general_numerical_approximation}
	Y_{0} = X_{0},
	\quad
	Y_{n+1} = Y_{n} + \Phi ( t_{n},Y_{n},h;\xi_{n} ),\quad n = 0,1,\ldots,N-1,
\end{equation}
where the $ \xi_{n} $ for $ n \geq 0 $ is independent of
$ Y_{_{0}}, Y_{_{1}},\ldots, Y_{_{n}} $,
$ \xi_{0}, \xi_{1},\ldots, \xi_{n-1} $.
Alternatively, one can write 
\begin{equation}
\label{eq:Yn-different-express}
	Y_{n+1} = Y ( t_n, Y_{n}; t_{n+1} )
	        = Y ( t_0, Y_{0}; t_{n+1} ),
	        \quad
	        Y_0 = X_0,
	        \quad n = 0,1,\ldots,N-1.
\end{equation}
{\color{black}
We aim to analyze the weak error of the one-step approximation in terms of
\begin{equation}
\big | \E \big [
  	 \varphi \big (
    	X(t_0,X_0;T)
  	    \big ) \big ]
  	     -\E \big[
  	      \varphi \big (
     	Y(t_0,Y_0; t_N)
  	   \big )\big ] \big |,
\end{equation}
where $\varphi \colon \R^d \rightarrow \R$ is called the test function.
%
%
To carry out the weak error analysis,
we introduce  the function $u \colon [0, T] \times \R^d \rightarrow \R$
  defined by
\begin{equation}
\label{eq:u-defn}
u(t,x)
= \E  \big[ \varphi \big(X(t,x;T)\big) \big],
\end{equation}
which is the unique classical solution of the associated Kolmogorov partial differential equation (PDE) under certain regularity assumptions on $f, g, \varphi$ (see \cite[Theorem 1.6.2]{Cerrai2001second} for more details).
Boundedness and differentiability properties of such
solutions $u$ to the Kolmogorov PDE play an essential role 
in the analysis of weak errors of numerical schemes.
%
In \cite{Milstein1985Weak} (see also \cite{Milstein2004stochastic}, \cite{Milstein2005Numerical}),
Milstein establishes a well-known weak convergence theorem
for general one-step approximation schemes of SDEs with sufficiently smooth, globally Lipschitz continuous coefficients and sufficiently smooth test functions with polynomial growth.
In \cite{Guyon2006euler}, the author discusses the classes of test functions $\varphi$ that can be allowed and also addresses how to interpret $\E \big[ \varphi \big(X(t,x;T)\big) \big]$ for such cases and has to assume strong conditions on $f$ and $g$. 
Another approach to prove weak convergence and, if necessary, reduce assumptions on $\varphi$ uses Malliavin calculus, see, e.g. \cite{Kohatsu2001weak}. A further viewpoint on the topic is provided by the concept of cubature on Wiener space, see [6], in particular for conditions on $\varphi$. Overall, the relations among the conditions on $\varphi$, the definition of weak convergence, and the existence and properties of $u$ are standard practices in the literature.
%
%

In this section, we revisit the general weak convergence theorem \cite{Milstein1985Weak,Milstein2004stochastic,
Milstein2005Numerical}.
Following the main idea of Milstein’s weak error analysis, 
we aim to reformulate a general weak convergence theorem for one-step approximations of SDEs
without global Lipschitz condition. \textit{First}, we assume on the one-step approximation error as follows.}

\begin{assumption}[One-step approximation error]
	\label{ass:u-local-error}
{\color{black}{For any $t \in [0, T]$ and $x \in \R^d$, $u$ given by \eqref{eq:u-defn} is well-defined with $|u( t, x) | < \infty$}. 
Additionally, there exists a Borel function $\nu \colon \R^d \rightarrow \R^+$ such that}
	\begin{equation} 
	\label{eq:u-local-error}
 \sup_{ t \in [0,T - h]}	\big|
	\E \big[ u(t + h , X(t, x; t+h) ) \big]
	- 
	\E \big[ u( t + h , Y(t, x; t+h) )  \big] 
	\big|
	\leq
	\nu (x) h^{p+1},
	\quad
	\forall x \in \R^d,
	\end{equation} 
where $h > 0$ and $ t + h \leq T$.
\end{assumption}
{\color{black}Under the globally Lipschitz condition, 
the Euler-Maruyama method satisfies \eqref{eq:u-local-error}
with $p=1$ and various higher-order weak schemes were proposed by \cite{Abdulle2012high,Abdulle2013weak,Komori2017weak,Milstein2004stochastic,Rossler2006rooted}, which  satisfy \eqref{eq:u-local-error} with $p \geq 2$.}

\textit{Second}, we need an assumption on bounded moments of numerical solutions.
\begin{assumption}[Bounded moments of numerical solutions]
\label{ass:numerical-moments}
The numerical approximations $\{Y_n\}_{0\leq n\leq N}$ 
uniquely determined by \eqref{eq:general_numerical_approximation}
is well defined in $\R^d$
and there exists a function
 $ \Upsilon \colon \R^d \rightarrow \R^+ $ 
 such that
\begin{equation}
\label{eq:asum_approximation_moment_bound-general}
\sup_{N \in \mathbb{N}}\sup_{0 \leq n \leq N }
\E [ | \nu ( Y_{n} ) | ] \leq \Upsilon (Y_0) < \infty.
\end{equation}
	Here the function $ \nu $ comes from 
	Assumption \ref{ass:u-local-error}.
\end{assumption}
{\color{black}
Following the idea of the well-known weak convergence
theorem originally proposed by Milstein (\cite[Theorem 2.1]{Milstein1985Weak}, \cite [Chapter 2]{Milstein2004stochastic}, \cite{Milstein2005Numerical}), we obtain a similar conclusion under Assumptions \ref{ass:SDEs-wellposed}, \ref{ass:u-local-error}, \ref{ass:numerical-moments}.
} {\color{black}
We present the proof in Appendix \ref{append:proof-fund}.}
\begin{thm}[A general  weak convergence theorem]
 \label{thm:fundamental_weak_convergence}
Let Assumptions \ref{ass:SDEs-wellposed},
\ref{ass:u-local-error} and
\ref{ass:numerical-moments}
hold.
Then the numerical solution produced by \eqref{eq:general_numerical_approximation} 
has a weak convergence rate of order $ p $ in the sense that
 \begin{equation}
   \label{eq:thm-convergence-rate}
  	\big | \E \big [
  	 \varphi \big (
    	X(t_0,X_0;T)
  	    \big ) \big ]
  	     -\E \big[
  	      \varphi \big (
     	Y(t_0,Y_0; t_N)
  	   \big )\big ] \big |
    	\leq  \color{black}\Upsilon (X_0) T h^{p},
  \end{equation}
where  the functions 
 $\varphi \colon 
   \R^d \rightarrow \R $ 
 and  $ \Upsilon \colon \R^d
        \rightarrow \R^+$ 
come from  Assumptions \ref{ass:u-local-error} and 
\ref{ass:numerical-moments}, respectively.	
\end{thm}

{\color{black}Before closing this section,
we give some remarks on the key assumption \eqref{eq:u-local-error}. For the weak error analysis, boundedness and differentiability properties of the function $u$ should be fully exploited.
Essentially all articles studying weak convergence of various numerical methods for SDEs 
discuss the above properties of $u$ as a solution of the associated Kolmogorov equation, even in the non-Lipschitz case. 
Examples of such articles can be found in the introduction.
%
%
} 

\section{A practical corollary on weak convergence of one-step approximation schemes}
\label{section:corollaries-weak-convergence}
\subsection{On Assumption \ref{ass:u-local-error} for one-step errors}
%
%
In this subsection, we reveal how Assumption \ref{ass:u-local-error} can be satisfied in the case that the test functions $\varphi$ obey certain smoothness and
polynomial growth conditions. {\color{black}Similar to 
Theorem 2.1 (b) in \cite{Milstein2004stochastic}, we first make the following assumptions on one-step errors.}
\begin{assumption}
	\label{ass:local-error}
	%
	There exist a measurable function $K \colon \R^d 
	\rightarrow [0, \infty)$  such that, for any $i_j \in \{1,2,\cdots, d \}$,
	\begin{align}
	\label{eq:error-of-general-one-step-approximation}
	\bigg| \E \Big[ \prod_{j=1}^s ( \delta_{X, x} )^{i_j} \Big]
	- \E \Big[ \prod_{j=1}^s ( \delta_{Y, x} )^{i_j} \Big] \bigg|
	\leq 
 {\color{black} K ( x )} h^{p+1},
	\,\,\, s=1,\ldots ,2p+1,
	\\
	\label{eq:estimate-of-general-one-step-approximation-solution}
	\bigg\| \prod_{j=1}^{2p+2}  ( \delta_{X, x} )^{i_j} \bigg\|_{L^2 (\Omega; \R)}
	\leq {\color{black} K ( x )} h^{p+1}, \quad
	\bigg\| \prod_{j=1}^{2p+2}  ( \delta_{Y, x} )^{i_j} \bigg\|_{L^2 (\Omega; \R)}
	\leq {\color{black} K ( x )} h^{p + 1},
	\end{align}
	where we denote
	\begin{equation} \label{eq:notation-DeltaX-DeltaY}
	( \delta_{X, x} )^{i_j} := X^{i_j}(t,x;t+h)-x^{i_j}, \quad
	( \delta_{Y, x} )^{i_j} := Y^{i_j}(t,x;t+h)-x^{i_j},
	\quad
	i_j \in \{1,2,\cdots, d \}.
	\end{equation}	
\end{assumption}
Next, we make further assumptions on the moment boundedness and mean-square differentiability properties of the solution to SDE \eqref{eq:SDE-X(t,x;T)}. 
\begin{assumption}
\label{ass:Xt-moment-Deritivate-moment}
Assume that the SDE \eqref{eq:SDE-X(t,x;T)} admits a unique  
	$ \{\mathcal{ F }_s\}_{ s \in [t, T] } $-adapted  solution $X(t, x; s ), 0 \leq  t \leq s  \leq T$
	such that for sufficiently large $ q_0  \geq 2 $
	\begin{equation}
	\label{eq:asum_solution_moment_bound}
			\sup_{ 0 \leq  t \leq s  \leq T }
			\E [| X(t, x; s )|^{q_0}] \leq C ( 1 + |x|^{ q_0 }).
		\end{equation}
Moreover, the solution $ X(t, x; T ) $  is $ 2p + 2 $ times mean-square differentiable with respect to the initial data $x$
and for sufficiently large $ q_1 \geq 2  $ and for any $ j = 1,\ldots,2p+2 $ it holds that
\begin{equation} \label{eq:ass-Derivative-X}
 \sup_
 { x \in \R^d }
 \E [ |D^j X (t,x; T)  |^{q_1} ]
 \leq C(T,q_1,j)
 < \infty,
\end{equation}
where $ C(T,q_1,j) $ is a constant depending on $ T ,q_1, j $.
\end{assumption}
{\color{black} We mention that Assumption \ref{ass:Xt-moment-Deritivate-moment} is important and essentially used to obtain boundedness and differentiability properties of the function $u (t , x)$ given by \eqref{eq:u-defn} (see the proof of Lemma \ref{lem:equivent_condition}). A natural question thus arises as to how this assumption can be satisfied in terms of conditions on coefficients $f,g$ of SDEs. In subsection \ref{subsec:conditions-coefficients}, we present explicit conditions on coefficients of SDEs to fulfill the Assumption \ref{ass:Xt-moment-Deritivate-moment}.}
For   simplicity, the letter
$ C $ is used to denote a generic positive
constant independent of time step size and
may vary for each appearance. Additionally, 
the moment boundedness of the numerical solution is needed
{\color{black}(see also Theorem 2.1 (d) in \cite{Milstein2004stochastic}).} 
\begin{assumption}
\label{ass:numerical-moments2}
The numerical approximations $\{Y_n\}_{0\leq n\leq N}$ 
produced by \eqref{eq:general_numerical_approximation}
is well-defined in $\R^d$
and for sufficiently large $ q_2  \geq 0 $, $\beta \geq 1$ it holds that
\begin{equation}
\label{eq:asum_approximation_moment_bound}
\sup_{N \in \mathbb{N}}\sup_{0 \leq n \leq N }
\E [| Y_{n}|^{q_2}] \leq C(1+|Y_0|^{ \beta q_2 }).
\end{equation}
\end{assumption}
{\color{black} Given \eqref{eq:asum_solution_moment_bound}, it would be natural to require that $\beta$ in \eqref{eq:asum_approximation_moment_bound} be equal to $1$. Indeed, \eqref{eq:asum_approximation_moment_bound} with $\beta = 1$ holds for some modified Euler schemes and the backward Euler method \eqref{eq:Problem_SDE} (see Theorem \ref{thm:bound_moment_modi_EM} and Proposition \ref{pro:bounded-moment-BEM}). However, this is not the case for some tamed-type methods or the balanced method (see subsection 2.1.4 in \cite{Tretyakov2013fundamental}).} 
With these assumptions as well as smoothness and
polynomial growth conditions on test functions $\varphi$, 
it is straightforward to verify 
that Assumption \ref{ass:u-local-error} holds.
\begin{lem}
\label{lem:equivent_condition}
{\color{black}{Let $\varphi \in \mathbb{H}^{2p+2}$ where, for some constants $L>0$ and $\chi \geq 1$,}}
\begin{equation}\label{eq:varphi-bound}
    |D^j \varphi (x)| \leq L (1 + |x|^{\chi} ),
    \quad
    j \in \{ 1, 2, ..., 2p+2\}.
\end{equation}
Let Assumptions \ref{ass:Xt-moment-Deritivate-moment},
\ref{ass:numerical-moments2}
be fulfilled with
$ q_0 = 2 \chi $, $ q_1 = 4 p + 4 $, 
$ q_2 = 2 \chi$
and let Assumption 
\ref{ass:local-error}
be satisfied.
Then $u$ defined by \eqref{eq:u-defn}  fulfills Assumption \ref{ass:u-local-error} 
with
$ \nu ( x ) = C ( 1 + | x |^ {\beta \chi} ) {\color{black} K ( x )}  $. 
\end{lem}
This lemma is proved in Appendix \ref{append:proof-one-step}.
{\color{black}The key idea originates from \cite{Milstein1985Weak,Milstein2004stochastic}.
Roughly speaking, one can use the chain rule to show that 
$u (t, x) = \E  \big[ \varphi \big(X(t,x;T)\big) \big] $ is $ 2p + 2 $ times differentiable with respect to $x$ and $u(t, \cdot ) \in \mathbb{H}^{2p+2}$ for any $t\in [0, T]$,
by noting that $ \varphi \in \mathbb{H}^{2p+2} $ is 
$ 2p + 2 $ times differentiable and 
the solution $X(t, x; s)$ of the SDE \eqref{eq:SDE-X(t,x;T)} is $ 2p + 2 $ times mean-square differentiable with respect to $x$. Then we employ the Taylor expansion of $u$ and rely on Assumption \ref{ass:local-error} to obtain the one-step error \eqref{eq:u-local-error} (consult Appendix \ref{append:proof-one-step} for more details).} 
Combining Theorem \ref{thm:fundamental_weak_convergence} with Lemma \ref{lem:equivent_condition} gives the following 
corollary on weak convergence rates of one-step approximation.
\begin{cor} 
 \label{co_local_error}
Let the assumptions of Lemma \ref{lem:equivent_condition}
hold and in addition let Assumption \ref{ass:numerical-moments} be satisfied with
$\nu ( x ) = C ( 1 + | x |^ { \beta \chi } ) {\color{black} K ( x )} $.
Then the numerical solution produced by \eqref{eq:general_numerical_approximation} 
has a global weak convergence rate of order $ p $, namely, 
 \begin{equation}
  	\big | \E \big [
  	 \varphi \big (
    	X(t_0,X_0;T)
  	    \big ) \big ]
  	     -\E \big[
  	      \varphi \big (
     	Y(t_0,Y_0; t_N)
  	   \big )\big ] \big |
    	\leq \Upsilon (X_0) T h^{p},
  \end{equation}
where the function $\Upsilon$ comes from Assumption \ref{ass:numerical-moments}.
\end{cor}
%
\subsection{Explicit conditions on coefficients of SDEs to fulfill Assumption \ref{ass:Xt-moment-Deritivate-moment}}
\label{subsec:conditions-coefficients}
In  Assumptions \ref{ass:Xt-moment-Deritivate-moment}, it is assumed that the underlying SDE \eqref{eq:SDE-X(t,x;T)} 
is well-posed with bounded moments and that the solution is mean-square differentiable with respect to the initial data.
%
%
{\color{black}It is a classical result that the solution $ X(t, x; s ) $ to the SDE \eqref{eq:SDE-X(t,x;T)}  is $ 2p + 2 $ times mean-square differentiable with respect to the initial data $x$ and satisfies all conditions of Assumption \ref{ass:Xt-moment-Deritivate-moment}, when the coefficients $f(x), g(x)$ of the SDE \eqref{eq:SDE-X(t,x;T)} are globally Lipschitz continuous and $2p+2$ times differentiable with respect to $x$.
}
{\color{black}To fulfill all conditions of Assumptions \ref{ass:Xt-moment-Deritivate-moment} in a non-globally Lipschitz setting,  we 
recall the framework of Cerrai \cite{Cerrai2001second} (see Hypotheses 1.1, 1.2 therein), where the drift and diffusion coefficients are assumed to be sufficiently smooth with polynomial growth and satisfy certain monotonicity conditions.} 
\begin{assumption}
\label{ass:coefficient_function_assumption}
\begin{enumerate}
    \item[(A1)]  
The drift coefficient function
$ f \in C^{2p+2} ( \mathbb{R}^{d}; \mathbb{R}^{d}) $
and there exists $ r \geq 0 $ such that
\begin{equation}
\label{ass:drift-f}
\sup_{x \in \mathbb{R}^{d}}
\frac{ | D^{\alpha} f(x) |}{1+ |x|^{2r+1-j}} < \infty , \quad | \alpha |=j,
\,\,\,  j \in \{ 0, 1, \ldots , 2p+2\};
\end{equation}
\item[(A2)] 
The diffusion coefficient function $ g \in C ^{2p+2} (  \mathbb{R}^{d};  \mathbb{R}^{d \times m} ) $
and there exists $ \rho \leq r $ such that
\begin{equation}
\label{ass:diff-g}
\sup_{x \in \mathbb{R}^{d}}
\frac{\| D^{\alpha} g (x) \|}{1 +
 | x |^{\rho-j}} < \infty ,
 \quad  | \alpha |=j,
 \,\,\,  j \in \{ 0, 1, \ldots , 2p+2\};
\end{equation}
\item[(A3)]
For all $\lambda_0>0$, there exists
$ C_{\lambda_0} \in \mathbb{R} $
such that
\begin{equation}
\label{eq:assum-f-g}
       \langle D  f(x) y , y  \rangle 
           + \lambda_0 \| D g(x) y\|^{2}
            \leq C_{ \lambda_0 } |y|^{2}, 
            \quad  
              \forall x,y \in \mathbb{R}^{d};
\end{equation}
\item[(A4)]
There exist $ a_1 > 0 $ and $ \gamma  , c_1 \geq 0 $ such that for any $ x , y \in \mathbb{R}^{d}$ it holds
\begin{equation}
\langle f(x+y)-f(x),y \rangle \leq
 - a_1 | y |^{2r+2} + c_1 (| x |^{\gamma} + 1).
\end{equation}
\end{enumerate}
\end{assumption}
\noindent
Examples of nonlinearity satisfying  Assumption \ref{ass:coefficient_function_assumption} can be found in 
\cite[Remark 1.1.1]{Cerrai2001second}.
In the proofs which follow we will need some implications of Assumption \ref{ass:coefficient_function_assumption}.
The Assumption (\textit{A1}) immediately implies
\begin{equation}\label{contronal:f}
  |f(x) - f(y) | \leq C (1 + |x|^{2r} + |y|^{2r})
  |x - y|, \quad \forall x,y \in \R^d,
\end{equation}
which gives the polynomial growth
\begin{equation}
\label{eq:f-growth}
  |f(x)| \leq C (1 + |x|^{2r+1}),
  \quad
  \forall  x \in \R^d . 
\end{equation}
Similarly, it follows from Assumption (\textit{A2}) that
\begin{equation}\label{contronal:g}
  \|g(x)- g(y)\| \leq
   C (1 + |x|^{\rho-1} + |y|^{\rho-1})
   |x - y|, \quad \forall x,y \in \R^d,
\end{equation}
and thus
\begin{equation}
  \|g(x) \| \leq C (1 + |x|^\rho),
  \quad
  \forall x \in \R^d .
\end{equation}
%
{\color{black}In addition, Assumption \ref{ass:coefficient_function_assumption} further ensures Khasminskii-type and one-sided Lipschitz property on $f$ and $g$, which are employed in strong approximation schemes for SDEs such as tamed and semi-implicit schemes (see, e.g., \cite{higham2002strong,Hutzenthaler2012Strong,Tretyakov2013fundamental} and references therein).}
More precisely, Assumption (\textit{A3})
implies one-sided Lipschitz property:
\begin{equation}
\label{eq:mono-condi}
\langle x - y , f(x) - f(y) \rangle
+ \lambda_0 \|g(x)-g(y)\|^{2} \leq C_{\lambda_0} | x - y |^{2},
\quad
\forall x,y \in \mathbb{R}^{d}.
\end{equation}
%
%
{\color{black}Under Assumption \ref{ass:coefficient_function_assumption}, Theorem 1.3.6 of \cite{Cerrai2001second} asserts that the SDE \eqref{eq:SDE-X(t,x;T)} is well-posed, with the solution is $2p+2$ times mean-square differentiable with respect to the initial data $x$.}

\begin{thm}{\color{black}\cite[Theorem 1.3.6]{Cerrai2001second}}
\label{thm:solution-mean-square-differentiable-uniform-bound}
Let Assumption \ref{ass:coefficient_function_assumption}  
be fulfilled.
Then  the SDE \eqref{eq:SDE-X(t,x;T)} admits a unique adapted
solution $ X(t, x; s), 0 \leq t \leq s \leq T $ in $ \R^d $, which is $ 2p+2 $ times mean-square differentiable
with respect to the initial data $ x \in \R^d$ and for any 
$ q \geq 1 $
and $ j = 1,\ldots,2p+2 $ it holds that
\begin{equation}
 \sup_
 { x \in \mathbb{R}^d }
 \E [ |D^j X ( t, x; s ) |^{q} ]
 \leq
 \phi_{j,q}(s-t)
 \leq C(T,q,j),
\end{equation}
for suitable continuous increasing functions
$ \phi_{j,q}(t) $, where $ C(T,q,j) $ is a 
constant that depends on $ T ,q,j $.
Moreover, for any $ q \geq 1 $ it holds that
\begin{equation}
 \label{eq:lem-q-moment-bound-sol}
 \sup_{ 0 \leq t \leq s \leq T } 
 \E [ | X(t, x; s) |^{q}
 ]
 \leq C(T,q)
  (1 + | x |^{ q }),
\end{equation}
and for any $ t, s_1, s_2 \in [0,T] $,
$ \gamma \geq 1 $, we have
\begin{equation} \label{eq:Xt-Xs-Holder}
  \| X ( t, x ; s_2 ) - X ( t, x; s_1 ) \|_{L^{\gamma}(\Omega;\mathbb{R}^d)}
 \leq C ( 1 + |x|^{2 r + 1} ) | s_2 - s_1 |^{\frac{1}{2}} .
\end{equation}
\end{thm}
%
{\color{black}Theorem \ref{thm:solution-mean-square-differentiable-uniform-bound} extends relevant results in the Lipschitz setting to a non-Lipschitz setting. Moreover,
this together with boundedness and differentiability properties of the test function $\varphi$ (i.e., $\varphi \in \mathbb{H}^{2p+2}$) helps establish smoothness estimates of $u$, i.e., $u (t, x) $ is $ 2p + 2 $ times differentiable with respect to $x$ and $u(t, \cdot ) \in \mathbb{H}^{2p+2}$ for any $t\in [0, T]$ (see the proof of Lemma \ref{lem:equivent_condition} for details). }
As a direct 
consequence of \eqref{contronal:f},
\eqref{contronal:g},
\eqref{eq:lem-q-moment-bound-sol}
and \eqref{eq:Xt-Xs-Holder}, we have, for any $ t, s_1, s_2 \in [0,T] $,
\begin{align}
\label{eq:fXt-fxs-Holder}
  \| f ( X ( t, x ; s_2 ) ) - f ( X ( t, x; s_1 ) ) \|_{L^{\gamma}(\Omega;\mathbb{R}^d)}
 & \leq C ( 1 + |x|^{4 r + 1} ) | s_2 - s_1 |^{\frac{1}{2}} 
 \\
 \label{eq:gXt-gxs-Holder}
 \| g ( X ( t, x ; s_2 ) ) - g ( X ( t, x; s_1 ) ) \|_{L^{\gamma}(\Omega;
 	      \mathbb{R}^{d \times m})} 
 & \leq C ( 1 + |x|^{2 r + \rho } ) | s_2 - s_1 |^{\frac{1}{2}} .
 \end{align}
Theorem \ref{thm:solution-mean-square-differentiable-uniform-bound} guarantees that all conditions in Assumption \ref{ass:Xt-moment-Deritivate-moment} hold,
for any $q_0, q_1 \geq 2$. As a direct consequence of
Corollary \ref{co_local_error}
and Theorem \ref{thm:solution-mean-square-differentiable-uniform-bound}, 
one can immediately get the following practical corollary on weak convergence of one-step approximation schemes.
{\color{black}
\begin{cor}[Weak convergence of one-step approximation under Assumption \ref{ass:coefficient_function_assumption}]
\label{cor-fundamental-weak-convergence}
Let $\varphi \in \mathbb{H}^{2p+2}$ 
and let Assumption 
\ref{ass:local-error}
be satisfied, where, for some constants $L>0$ and $\varkappa \geq 0$,
\begin{equation}
K ( x ) \leq  L ( 1 + |x|^{\varkappa} ).
\end{equation}
Let Assumption \ref{ass:numerical-moments2}
be fulfilled with  $ q_2 = (2 \chi) \vee ( \beta \chi +  \varkappa) $ and let Assumption
\ref{ass:coefficient_function_assumption} hold.
Then the numerical solution produced by \eqref{eq:general_numerical_approximation}
has a global weak convergence rate of order $ p $, namely, 
 \begin{equation}
  	\big | \E \big [
  	 \varphi \big (
    	X(t_0,X_0;T)
  	    \big ) \big ]
  	     -\E \big[
  	      \varphi \big (
     	Y(t_0,Y_0; t_N)
  	   \big )\big ] \big |
    	\leq C 
     \big( 1 + | X_0 |^{ \beta ( \beta \chi + \varkappa ) }  \big) h^{p}.
  \end{equation}
\end{cor}
}
\section{One-step weak errors for the Euler approximation}
\label{sec:one-step-euler}
To analyze weak convergence rates of numerical schemes
considered in this work, we introduce the following one-step Euler approximation:
\begin{equation}
\label{eq:one_step_Euler_Maruyama_method}
 Y_E(t,x;t+h)
      = x +
       \int_t^{t+h}f(x) 
       \, \dd s
        + \int_t^{t+h} g(x) 
        \, \dd W(s)
\end{equation}
and
denote
\begin{equation}
 (\delta_{Y_E,x})^{i_j} := Y^{i_j}_E(t,x;t+h)-x^{i_j},
 \quad
  i_j \in \{ 1,2,\ldots,d \},
\end{equation}
where by $ x^{i_j} $ we mean the $ i_j $-th coordinate of the vector $ x $.
{\color{black}
The weak error is decomposed into two parts:
one being the one-step weak error of the Euler approximation 
\eqref{eq:one_step_Euler_Maruyama_method} (Lemma \ref{lem:error-one-step-Euler-method}) and the other  
being the discrepancy between the Euler approximation and the considered schemes (see Lemma \ref{modified_one_error} 
and Lemma \ref{FBDE_one_step_estimation} for more details).}
We will need the following lemma to show the one-step weak weak error of 
\eqref{eq:one_step_Euler_Maruyama_method}. 

\begin{lem}[Moment bounds for one-step approximations]
\label{lem:norm-estimate-Euler-method}
Let Assumptions 
(\textit{A1}) and (\textit{A2}) be satisfied with $ p = 1$. For any $ q \geq 1$, we have
\begin{align}
\label{eq:norm-error-one-step-Euler-method}
\|
\delta_{X,x} -  \delta_{Y_E,x}
\|_{L^{2q}(\Omega;\mathbb{R}^d)}
& \leq C(1 + |x|^{4r + 1})h,
\\
\label{eq:norm-estimate-one-step-Euler-method}
\|
\delta_{Y_E,x}
\|_{L^{2q}(\Omega;\mathbb{R}^d)}
& \leq C(1+|x|^{2r+1})h^{\frac{1}{2}},
\\
\label{eq:norm-estimate-one-step-exact-solution}
\|
\delta_{X,x}
\|_{L^{2q}(\Omega;\mathbb{R}^d)}
& \leq C(1+|x|^{2r+1})h^{\frac{1}{2}}.
\end{align}
Here $r$ is defined in \eqref{ass:drift-f}
\end{lem} 

Lemma \ref{lem:norm-estimate-Euler-method} leads to the desired one-step weak errors for the Euler approximation \eqref{eq:one_step_Euler_Maruyama_method}.

\begin{lem}[One-step weak error of the Euler approximation]
\label{lem:error-one-step-Euler-method}
Under the same condition of Lemma \ref{lem:norm-estimate-Euler-method}.
Then for $i_j \in \{1,2,\cdots, d \}$ we have
\begin{equation}
\label{eq:error-one-step-Euler-method}
   \bigg| \E \Big [ 
     \prod_{j=1}^s (\delta_{X,x})^{i_j}
       - \prod_{j=1}^s (\delta_{Y_E,x})^{i_j} 
         \Big] \bigg|
           \leq
            C (1+|x|^{ 8r+3 })h^2 , \,\, s=1,2,3.
\end{equation}
\end{lem}
The proofs of these two lemmas can be found in Appendix \ref{appendix:proof_EM_one_step_error}.
%

\section{Weak convergence rates of modified Euler approximations}
\label{sec:weak-rate-MEM}
In this section, we prove weak convergence rates of several modified Euler schemes. By Corollary \ref{cor-fundamental-weak-convergence}, we need to prove boundedness of moments and analyze the local error. 
Given a uniform mesh on $ [0,T] $ 
with $ h = \frac{T}{N} > 0 $ 
equispaced time step sizes, the following numerical scheme is defined:
\begin{equation}
\label{eq:modified_EM_scheme}
Y_{n+1}
  =Y_{n} + \bar{f}_{h}(Y_n)h
   +\bar{g}_{h}(Y_n) \Delta W_{n},
\quad
Y_{0} = X_{0},
\end{equation}
where
 $ \Delta W_{n}:= W(t_{n+1})- W(t_{n}), $ 
 $ n \in \{ 0,1,2,\ldots,N-1 \} $ and 
 $\bar{f}_{h} \colon \Omega \times \R^d
\rightarrow \R^d,$
$\bar{g}_{h} \colon \Omega \times \R^d
\rightarrow \R^{d \times m}$ are some measurable functions.
 The weak convergence rate of this scheme will be analyzed in general.
 To apply the previously obtained corollary on weak convergence, the first step is to discuss the boundedness of the moments of the numerical approximations \eqref{eq:modified_EM_scheme}.
\subsection{Boundedness of moments of modified Euler approximations}
 In this subsection we discuss moment boundedness for numerical solutions of the modified Euler scheme \eqref{eq:modified_EM_scheme}.  
 As noted in \cite{Tretyakov2013fundamental}, the bounded moments for numerical solutions are usually discussed case by case.
 Next, we illustrate the moment boundedness using the conditions from \cite[Lemma 2]{Sabanis2016euler}, 
 which require $\bar{f}_{h}$ and $\bar{g}_{h}$ to
 be deterministic and satisfy the following properties.  
 {\color{black}
 \begin{assumption}
   \label{con:modified_EM_dri_diff}
Let $\bar{f}_{h} \colon \R^d
\rightarrow \R^d,$
$\bar{g}_{h} \colon \R^d
\rightarrow \R^{d \times m}$ be deterministic functions satisfying the following properties.
\\
(1)
 There exists $\theta \in (0,\frac{1}{2}]$ such that, for $\forall x \in \R^d,$
  \begin{align}
\label{con:modified_EM_sublinear_fhx}
       |\bar{f}_{h}(x) |  
      & \leq
        \min \big( Ch^{-\theta}(1+|x| ), |{f(x)} | \big),  \\
\label{con:modified_EM_suplinear_ghx}
      \| \bar{g}_{h}(x) \|^2
       & \leq
          \min \big(  Ch^{-\theta}(1+|x|^2),
          \|g(x)\|^2).
\end{align}
(2)
There exists a sufficiently large constant $ \mu > 0 $ 
such that 
  \begin{equation}
      \label{con:modified_EM_onesided_lip}
      \langle x, \bar{f}_{h}(x) \rangle 
      + \mu \| \bar{g}_{h}(x) \|^2 
      \leq
      C(1+|x|^2), \quad x \in \R^d.
  \end{equation}
 \end{assumption}
 }
The moment bounds of the numerical approximations produced by \eqref{eq:modified_EM_scheme} are proved under Assumption 
\ref{con:modified_EM_dri_diff} (see \cite[Lemma 2]{Sabanis2016euler}). 
 Specifically, the author obtained the moment bounds of the following continuous version  of the numerical scheme 
\eqref{eq:modified_EM_scheme} 
\begin{equation}
	\label{eq:numerical_solution_continuous}
	Y_t = Y_0 + \int_{0}^{t}
    	\bar{f}_{h}(Y_{	\lfloor s \rfloor_h } ) 
	    \, \dd s
		+ \int_{0}^{t} 
	\bar{g}_{h}(Y_{	\lfloor s \rfloor_h } ) 
	  \, \dd W(s), 
	  \quad
	  t \in [0,T].
\end{equation}
where  $
	\lfloor t \rfloor_h
	   := t_{i}
    $
	for 
     $ t 
	\in [t_{i}, t_{i+1} ), 
	\, i \in\{ 0, 1, \ldots, N-1 \}$. 
\begin{thm}[\cite{Sabanis2016euler}]
 \label{thm:bound_moment_modi_EM}
Let Assumption (\textit{A1}) be satisfied with $ p = 1$ and in addition assume that the Assumption (\textit{A3}) and Assumption \ref{con:modified_EM_dri_diff}
hold. 
 Then for any $ q \geq 2 $, the approximation process
 $ Y_{ t } $
produced by 
\eqref{eq:numerical_solution_continuous}
obeys
\begin{equation}
      \sup_{ t \in [0,T] }
    \E [ | Y_{t}|^{q} ]
       \leq
         C ( 1+| X_0 |^ {q} ).
\end{equation}
Here $C$ depends on $T$, $q$
and is independent of $h$.
\end{thm}
 Even though the conditions in Assumption \ref{con:modified_EM_dri_diff} are general, 
 there are many modified Euler type schemes that do not satisfy these conditions, see, e.g., the tamed scheme \eqref{eq:numerical-approximation-DTE} and
 the balanced type schemes \eqref{eq:numerical-approximation-BS}, \eqref{eq:numerical-approximation-BTS}.
 We mention that Assumption \ref{con:modified_EM_dri_diff} is  only  used to get moment bounds of the numerical approximations. Once the moment  bounds of numerical approximations are available, we do not require the considered schemes to satisfy Assumption \ref{con:modified_EM_dri_diff}.  
  In Section \ref{ssec:example-mod-euler}, we present some examples of modified Euler schemes and apply our weak convergence theorem 
 to obtain their weak convergence rates. 

\noindent

\subsection{Weak convergence rates of the modified Euler approximations}
According to Corollary \ref{cor-fundamental-weak-convergence},
we still need to verify Assumption \ref{ass:local-error}, with the one-step weak convergence rate in
\eqref{eq:error-of-general-one-step-approximation} and \eqref{eq:estimate-of-general-one-step-approximation-solution}. To this aim, let us consider the one-step approximation of  the numerical method \eqref{eq:modified_EM_scheme},
given by
\begin{equation}
\label{eq:one_step_modified__Euler_method}
Y(t,x;t+h)
    = x +
     \bar{f}_{h}(x)
     h
      +
      \bar{g}_{h}(x)
      \big( W( t+ h) - W ( t ) \big).
\end{equation}
For the purpose of the weak error analysis for the modified Euler schemes, we make the following assumptions, describing the discrepancy between $f, g$ and $\bar{f}_{h}, \bar{g}_{h}$.
{\color{black}
\begin{assumption}
\label{ass:error-drift-diff-modi-EM}
For $\bar{f}_{h} \colon \Omega \times \R^d
\rightarrow \R^d,$
$\bar{g}_{h} \colon \Omega \times \R^d
\rightarrow \R^{d \times m}$,
 there exist $ \mathsf{a}, \mathsf{b},
 r_1, \rho_1 \geq 0 $ such that,
 for any $ q \geq 1,$
  \begin{align}
     \label{con:modified_EM_sup_lin_f}
      \| \bar{f}_{h}(x) - f(x)
       \|_{L^{2q}(\Omega;\mathbb{R}^d)}
      & \leq
         Ch^{\mathsf{a}}(1+|x|^{r_1}), 
             \quad
        x \in  \mathbb{R}^d, 
\\
     \label{con:modified_EM_sup_lin_g}
      \| \bar{g}_{h}(x) - g(x)
       \|_{L^{2q}(\Omega;\R^{d \times m})}
       & \leq
            Ch^{\mathsf{b}} ( 1 + |x|^{\rho_1} ),
             \quad
           x \in  \mathbb{R}^d,
\\
\label{con:modified_EM_sup_lin_g_expectation}
\big| \E \big[ \big( 
        g(x)
        - \bar{g}_{h}(x)
             \big) \big( W( t + h) - W ( t ) \big)
    \big]
\big|
    & \leq
    Ch^{\mathsf{b} + 1} ( 1 + |x|^{\rho_1} ),
    \quad
    x \in \mathbb{R}^d, h > 0.
\end{align}
\end{assumption}
}
{\color{black}Assumption \ref{ass:error-drift-diff-modi-EM} indicates that different choices of $\bar{f}_{h}, \bar{g}_{h}$ would lead to different values of parameters $\mathsf{a},\mathsf{b}$, which determine the weak convergence
rates of the schemes (see Theorem \ref{thm:MEM-global-conver-rate} and subsection \ref{ssec:example-mod-euler}).
} 
Under this assumption, we can get the strong error estimates for the one-step approximations, which are needed for the desired weak error estimates.
\begin{lem}
   \label{lem:norm-esti-one-step-modified-EM}
Under the Assumptions (\textit{A1}), (\textit{A2})
and Assumption \ref{ass:error-drift-diff-modi-EM}, then  
for any $ q \geq 1 $ it holds that
\begin{align}
  \label{eq:norm-esti-one-step-EM-MEM}
    \|
          \delta_{Y_E,x} -  \delta_{Y,x} 
     \|_{L^{2q}(\Omega;\mathbb{R}^d)}
     & \leq
       C (1+|x|^{r_1}) h^{\mathsf{a}+1}
  + 
    C (1+|x|^{\rho_1}) h^{\mathsf{b}+\frac{1}{2}},
\\
  \label{eq:norm-esti-one-step-TEM}
      \|
           \delta_{Y,x}
       \|_{L^{2q}(\Omega;\mathbb{R}^d)}
    & \leq 
           C  
 (1+|x|^{2r+r_1+\rho_1+1})h^{\frac{1}{2}}.
\end{align}
\end{lem}
The detailed proof of Lemma \ref{lem:norm-esti-one-step-modified-EM} is presented in the Appendix \ref{appendix:moment_esti_one_step_EM_MES}.
With the aid of Lemma \ref{lem:norm-esti-one-step-modified-EM}, the next lemma gives the local error estimates of the modified Euler scheme 
\eqref{eq:modified_EM_scheme}.
\begin{lem}[One-step weak error estimates]
\label{modified_one_error}
Under the same conditions of Lemma \ref{lem:norm-esti-one-step-modified-EM}, 
then the one-step approximation 
\eqref{eq:one_step_modified__Euler_method}
of the numerical method \eqref{eq:modified_EM_scheme} satisfies
\begin{align}
\label{sheme:modified_EM_general_one_step_error}
       \bigg| \E \Big [ 
         \prod_{j=1}^s (\delta_{X,x})^{i_j}
       - 
        \prod_{j=1}^s (\delta_{Y,x})^{i_j} 
          \Big] \bigg|
        & \leq
         C (1+ |x|^{8r+3\rho_1+3r_1+3}) h^{2 \land (\mathsf{a}+1) \land (\mathsf{b}+1)},
          \,\, s=1,2,3,
       \\
\label{sheme:solution_one-step approximation2}
   \bigg\| \prod_{j=1}^{4}  
     ( \delta_{X, x} )^{i_j} 
   \bigg\|_{L^2 (\Omega; \R)}
& \leq  C (1+|x|^{8r+4}) h^2,
       \\
\label{eq:s4_norm_esti_one_step_modified_EM}
  	\bigg\| \prod_{j=1}^{4}  
  	      ( \delta_{Y, x} )^{i_j} 
  	\bigg\|_{L^2 (\Omega; \R)}
& \leq  C (1+|x|^{ 8r+4r_1+4\rho_1+4 })  h^2.
\end{align}
\end{lem}
The proof of this lemma is given in Appendix \ref{appendix:moment_esti_one_step_MES_sol}.
By Theorem \ref{thm:bound_moment_modi_EM}, Lemma \ref{modified_one_error}, 
and Corollary \ref{cor-fundamental-weak-convergence}, we are able to identify the weak convergence rates of the modified numerical methods \eqref{eq:modified_EM_scheme} as follows.
\begin{thm}[Weak convergence rate for the method \eqref{eq:modified_EM_scheme}]
 \label{thm:MEM-global-conver-rate}
 Let the Assumption \ref{ass:coefficient_function_assumption} hold with $p=1$ and in addition assume that the Assumptions \ref{con:modified_EM_dri_diff}, 
 \ref{ass:error-drift-diff-modi-EM} hold.
Then the approximation method \eqref{eq:modified_EM_scheme} has a global weak convergence rate of order $(1 \land \mathsf{a} \land \mathsf{b}) $: for $ \varphi \in \mathbb{H}^{4} $, 
%
\begin{equation}
  \label{eq:tam conver results}
	\big |\E \big [
	 \varphi \big (
	  X(t_0,X_0;T)
	   \big ) \big ]
	  -\E \big[
	 \varphi \big (
	  Y( t_0,Y_0;T )
	  \big )\big ] \big |
	 \leq 
	   C 
	   \big( 1+ |X_0|^{ (8r+4\rho_1+4r_1+\chi+4) }
	    \big) h^{1 \land \mathsf{a} \land \mathsf{b}},
 \end{equation}
 where  $r$ and $\chi$ are defined in \eqref{ass:drift-f} and \eqref{eq:varphi-bound}, respectively.
If all conditions above are satisfied except that Assumption \ref{con:modified_EM_dri_diff} is replaced
by boundedness of moments of numerical solutions in Assumption \ref{ass:numerical-moments2}, then
the approximation method \eqref{eq:modified_EM_scheme} has a global weak convergence rate of order $(1 \land \mathsf{a} \land \mathsf{b}) $.
\end{thm}
\subsection
{Examples of   modified Euler schemes}\label{ssec:example-mod-euler}
 Here are some examples of modified Euler type schemes 
 \eqref{eq:modified_EM_scheme}, where $\bar{f}_h$ and $\bar{g}_h$ are explicitly given.
 Furthermore, Theorem \ref{thm:MEM-global-conver-rate}
 is applied to obtain the expected weak convergence rates
 of the considered schemes. 
\begin{exm}(Fully tamed Euler (FTE) \cite{Sabanis2016euler}) 
  \label{example-sabanis-model1}
The author of \cite{Sabanis2016euler} proposes two drift-diffusion fully tamed Euler schemes \eqref{eq:modified_EM_scheme}, where $\bar{f}_{h}, \bar{g}_{h}$ are given as follows.
\par
\vspace{0.1cm}
Model 1 (FTE1):
 \begin{equation}
     \label{eq:numerical-approximation-FTE1}
      \bar{f}_{h}(x)  =  
      \frac{f(x)}{1+h^{\alpha_1}|f(x)|+ h^{\alpha_2} \|g(x)\|^2},  \quad 
      \bar{g}_{h}(x)  =  \frac{ g(x)}  {1+h^{\alpha_1}|f(x)|+ h^{\alpha_2} \|g(x)\|^2}, \quad
      0<\alpha_1,\alpha_2\leq \tfrac{1}{2}.
 \end{equation}
One immediately observes that Assumption \ref{con:modified_EM_dri_diff} is satisfied with $ \theta = \alpha_1 \land \alpha_2 \leq \tfrac12. $
%
Also, by Assumption (\textit{A1}), (\textit{A2}), one can easily check that Assumption 
\ref{ass:error-drift-diff-modi-EM} is satisfied with 
\[r_1 = 4r+2, \, \rho_1 =2r+\rho+1, \, \mathsf{a} = \mathsf{b}=\alpha_1 \land \alpha_2.\]
From Theorem \ref{thm:MEM-global-conver-rate}, it follows that the weak convergence rate is $\alpha_1 \land \alpha_2 \leq \tfrac12$. We also present numerical results  (Examples \ref{ex:numerical_ex_nonlinear_diff} and \ref{ex:numerical_ex_linear_diff}) and verify the half-order convergence by taking $ \alpha_1 = \alpha_2 = \frac{1}{2}.$
\par
\vspace{0.1cm}
Model 2 (FTE2):
 \begin{equation}
  \label{eq:numerical-approximation-FTE2}
       \bar{f}_{h}(x)  =  
      \frac{f(x)}{1+h^{\vartheta} |x|^{2r} }, \quad 
      \bar{g}_{h}(x)   =  \frac{ g(x)}  {1+h^{\vartheta}|x|^{2r} },
     \quad 0< \vartheta \leq \tfrac{1}{2}.
  \end{equation}
Evidently, Assumption \ref{con:modified_EM_dri_diff} is fulfilled with $ \theta = \vartheta \leq \tfrac12. $  
%
In view of Assumptions (\textit{A1}), (\textit{A2}), 
we observe that Assumption 
\ref{ass:error-drift-diff-modi-EM} is satisfied with 
\[r_1 = 4r+1, \, \rho_1 =2r+\rho, \, \mathsf{a} = \mathsf{b}= \vartheta.\]
By Theorem \ref{thm:MEM-global-conver-rate}, the weak convergence rate is $\vartheta \leq \tfrac12.$ We also present numerical results (Examples \ref{ex:numerical_ex_nonlinear_diff} and \ref{ex:numerical_ex_linear_diff}) and verify the half-order convergence by taking $ \vartheta = \frac{1}{2}.$
\end{exm} 
 \begin{exm}[Modified Euler scheme (MES)]
  \label{example-new-tem}
In this work we propose a new modified Euler scheme, where
   \begin{equation}
      \label{eq:numerical-approximation-MES}
      \bar{f}_{h}(x) = \frac{f(x)}
                   {1+h|f(x)|^2},\quad 
       \bar{g}_{h}(x) = \frac{g(x)}
                   {1+h|f(x)|^2}.
    \end{equation}
%
\color{black}{In a similar way, one then observes that the Assumption \ref{con:modified_EM_dri_diff} is satisfied with $\theta=\frac{1}{2}$ and furthermore that due to Assumption (\textit{A1}) and (\textit{A2}), the Assumption 
\ref{ass:error-drift-diff-modi-EM} with 
\[r_1 = 6r+3, \, \rho_1 = 4r+\rho+2, \, \mathsf{a}= \mathsf{b}=1.\]}
 %
According to Theorem \ref{thm:MEM-global-conver-rate}, we derive the weak convergence rate of order $1$. We also present numerical results (Examples \ref{ex:numerical_ex_nonlinear_diff} and \ref{ex:numerical_ex_linear_diff}) to verify the first-order weak convergence.
 \end{exm}

Below, we discuss schemes which does not satisfy the Assumption \ref{con:modified_EM_dri_diff} but the moment boundedeness of numerical solutions are available in the literature
 \cite{Hutzenthaler2012Strong,Zhang2017preserving,Tretyakov2013fundamental}.
 \begin{exm}[Drift-tamed Euler (DTE) \cite{Hutzenthaler2012Strong}]
 In \cite{Hutzenthaler2012Strong},  the diffusion coefficient $g$ is assumed to be globally Lipschitz continuous. 
 In this setting, the diffusion coefficient $g$ does not need to be tamed and a drift-tamed Euler scheme is introduced, where
 \begin{equation}
      \label{eq:numerical-approximation-DTE}
      \bar{f}_{h}(x) = \frac{f(x)}
                   {1+h|f(x)| },\quad 
       \bar{g}_{h}(x) = g(x). 
    \end{equation}
The boundedeness of moments of the numerical approximations
can be found in \cite[Lemma 3.9]{Hutzenthaler2012Strong}.
Using similar arguments as above, one can show Assumption \ref{ass:error-drift-diff-modi-EM} is satisfied with 
$ r_1 = 4r+2, \mathsf{a} = 1.$
By Theorem \ref{thm:MEM-global-conver-rate}, we get the weak convergence rate of order $1$, which is confirmed by numerical experiment (Example \ref{ex:numerical_ex_linear_diff}).
 \end{exm}
 \begin{exm} \label{ex:balanced_scheme}
 [Balanced scheme (BS) \cite{Zhang2017preserving}]
 The authors of \cite{Zhang2017preserving} propose an explicit balanced scheme, where
 \begin{equation}
      \label{eq:numerical-approximation-BS}
      \bar{f}_{h}(x) = h^{-1} \tanh(h f(x)),\quad      
       \bar{g}_{h}(x) =  h^{-1/2} \tanh(h^{1/2} g(x)).
    \end{equation}
 As shown by \cite[Lemma 3.1]{Zhang2017preserving}, the numerical solutions have bounded moments.
By Taylor's formula and Assumptions (\textit{A1}), (\textit{A2}),  one can readily verify 
that Assumption \ref{ass:error-drift-diff-modi-EM} is satisfied with $ r_1 = 4r+2 $, $ \rho_1 = 3\rho $ and $\mathsf{a}= \mathsf{b}=1$.
According to Theorem \ref{thm:MEM-global-conver-rate}, the weak convergence rate is one. We also present numerical results (Examples \ref{ex:numerical_ex_nonlinear_diff} and \ref{ex:numerical_ex_linear_diff}) and verify the first-order convergence.
\end{exm} 
\begin{exm}
 \label{ex:balanced_type_scheme}
 [Balanced-type scheme (BTS) \cite{Tretyakov2013fundamental}]
The authors consider 
 a balanced-type scheme, where
 \begin{equation}
      \label{eq:numerical-approximation-BTS}
      \bar{f}_{h}(x)  =  
      \frac{f(x)}{1+h|f(x)|+ |g(x) [ W (t+h) - W (t) ] |},  
      \:
      \bar{g}_{h}(x)  =  \frac{g(x)}  {1+h|f(x)|+ |g(x) [ W (t+h) - W (t) ] |}.
    \end{equation}
Using similar arguments shows that Assumption \ref{ass:error-drift-diff-modi-EM} is satisfied with 
$ r_1 = 4r+2, \rho_1 = 2r+\rho+1,
\mathsf{a} = \mathsf{b} = \frac{1}{2}.$
By Theorem \ref{thm:MEM-global-conver-rate}, we obtain the weak convergence rate of order $\frac{1}{2}$, which is confirmed by numerical result (Example \ref{ex:numerical_ex_nonlinear_diff}). 
\end{exm} 
\section{Weak convergence rate of the backward Euler method}
\label{sec:weak-rate-BEM}
As another application of  the corollary on weak convergence,
we analyze the weak convergence rate of the well-known backward Euler method (BEM) applied to SDEs \eqref{eq:Problem_SDE}:
\begin{equation}
\label{eq:numerical-approximation-BEM}
        Y_{n+1}
      =
        Y_{n} + f(Y_{n+1})h + g(Y_{n})\Delta W_{n},
        \quad
        Y_0 = X_0,
\end{equation}
where 
    $
       \Delta W_{n}
     :=
        W(t_{n+1})- W(t_{n})$, 
       $ n \in \{ 0,1,2,\ldots,N-1 \}. $
Then the one-step approximation of
\eqref{eq:numerical-approximation-BEM} is 
\begin{equation}
\label{eq:numerical-appro-one-step-BEM}
         Y(t,x;t+h) 
       = 
         x 
       + \int_{t}^{t+h} f (Y(t,x;t+h)) 
       \, \dd s
       + \int_{t}^{t+h}g(x) 
       \, \dd W(s).
\end{equation}
Denote 
$ \delta_{X,x}:= X(t,x;t+h)-x, \,
  \delta_{Y,x}:= Y(t,x;t+h)-x.$
Then by 
\eqref{eq:one_step_Euler_Maruyama_method}, 
\begin{equation}
\begin{split}
\label{eq:Rf-defn}
      Y(t,x;t+h)
   =
       Y_E(t,x;t+h)+ R_f ,
       \quad
R_f :=  \int_t^{t+h} f(Y(t,x;t+h))-f(x) \, \dd s.
\end{split}
\end{equation}
{\color{black}
Before proceeding further, we shall establish the well-posedness of the implicit scheme \eqref{eq:numerical-approximation-BEM}.
\begin{prop}
Let the condition (\textit{A3}) in 
Assumption \ref{ass:coefficient_function_assumption} be satisfied and let $h C_{\lambda_0} < 1 $. Then the implicit scheme \eqref{eq:numerical-approximation-BEM} 
has a unique adapted solution $\{ Y_n \}_{n = 1}^N$ 
with probability one.
\end{prop}
The proof can be found in \cite[Theorem 3.1]{andersson2017mean}.
}
As already established in 
\cite[Corollary 2.27]{Hutzenthaler2015Numerical},
$q$-th moments of the BEM are uniformly bounded under Assumption (\textit{A3}).
\begin{prop}
 \label{pro:bounded-moment-BEM}
Suppose Assumption (\textit{A3}) holds. Then there exists 
   $ C > 0,$ 
independent of h such that for all 
   $ N \in \mathbb{N} $
and 
   $ n=0,1,\ldots,N, $
\begin{equation}
  \label{eq:bounded-moment-BEM}
    \sup_{N \in \mathbb{N}} \sup_{0 \leq n \leq N }
         \E [ | Y_{n}|^{q} ] 
           \leq
             C \big(1+ |X_0|^{q} \big),
        \,\,\,
        \forall \,
          q
         \geq 0,
\end{equation}
where
    $ \{Y_{n}\}_{0 \leq n \leq N } $ 
is given in  
\eqref{eq:numerical-approximation-BEM}.
\end{prop}
Equipped with the moment bound, we continue to analyze the weak convergence rate of the BEM \eqref{eq:numerical-approximation-BEM}
by means of Corollary \ref{cor-fundamental-weak-convergence}.
{\color{black}
The local weak errors for the BEM is similar to that for the modified Euler schemes. At first, we prove the following lemma on one-step strong error bounds of the BEM.
}
\begin{lem}
\label{eq:norm-esti-one-step-BEM}
Under the same condition of Proposition \ref{pro:bounded-moment-BEM}, letting Assumptions (\textit{A1}) and (\textit{A2}) hold with $p=1,$ then for any $ q \geq 1 $ the following inequalities hold
for the backward Euler scheme
\eqref{eq:numerical-approximation-BEM}:
\begin{align}
    \label{eq:norm-one-step-BEM}
     	\|
     	\delta_{Y,x}
    	\|_{L^{2q} (\Omega; \R^d) }
   & \leq 
      C (1+|x|^{2r+1})h^{\frac{1}{2}},
          \\
    \label{eq:norm-one-step-Rf}
          \| R_f  \|_{L^{2q} (\Omega; \R^d) }
  &  \leq
        C(1+|x|^{4r+1})h^{\frac{3}{2}},
         \\
     \label{eq:norm-error-one-step-solution-BEM}
           \|  \delta_{X,x}
        -
             \delta_{Y,x}
            \|_{L^{2q} (\Omega; \R^d) }
     & \leq 
             C (1+|x|^{4r+1})h,
\end{align}
where $R_f$ is defined in \eqref{eq:Rf-defn}.
\end{lem}
\textit{Proof of Lemma 
 \ref{eq:norm-esti-one-step-BEM}.} 
Repeating the same lines 
  in the proof of
   \eqref{eq:error_norm_one_step_EM}
   and
   \eqref{eq:bounded-moment-BEM}
results in
\begin{eqnarray}
  \label{eq:esti_one_step_norm_BEM}
   \E \big[| (\delta_{Y,x})|^{2q} \big]     
  & = &
     \E \Big[ \Big|
      \int_t^{t+h} 
       	f \big( Y(t,x;t+h) \big)
        \, \dd s  
   + 
        \int_t^{t+h} g(x) 
         \, \dd W(s)    \Big|^{2q} \Big]    \notag  \\
 & \leq &
       C h^{2q-1}
        \int_t^{t+h}  
          \big\| f \big( Y(t,x;t+h) \big)
          \big\|_{L^{2q} (\Omega, \R^d) }^{2q}
           \, \dd s 
  +  
     C h^{q-1} \int_t^{t+h} 
           	\big\| g(x)
            \big\|_{L^{2q} (\Omega, \R^{d \times m}) }
             ^{2q}
            \, \dd s     \notag \\
 & \leq &
       C (1+|x|^{ 2q(2r+1) }) h^{q},
\end{eqnarray}
which confirms
\eqref{eq:norm-one-step-BEM}.
To prove \eqref{eq:norm-one-step-Rf},
we use the H\"older inequality,
\eqref{eq:fXt-fxs-Holder}
and
 \eqref{eq:norm-one-step-BEM}
to derive that
\begin{eqnarray}
    \E \big[ \big|
    	R_f 
    	\big|^{2q} \big]   
     & \leq &
       C h^{2q-1}
       \int_t^{t+h}
        \big \|f(Y(t,x;t+h))- f(x)
        \big \|_{L^{2q} (\Omega, \R^d) }
          ^{2q}
         \, \dd s      \notag \\
    & \leq &
         C (1+|x|^{2q(4r+1)}) h^{3q}. 
\end{eqnarray}
Also, one can follow the same lines in the proof
of \eqref{eq:norm-error-one-step-Euler-method} 
to show that
\begin{eqnarray}
 \label{es:error_norm_one_step_solu_BEM}
    & &
        \E 
        \big[|
            \delta_{X,x}-\delta_{Y,x}|^{2q} 
          \big]      \notag \\
  & = &
       \E \Big[  \Big | 
           \int_t^{t+h} f \big( X(s) \big)
          - 
            f \big(Y(t,x;t+h)\big) 
            \, \dd s   
         + 
            \int_t^{t+h} g \big(X(s) \big)
         - 
             g(x)
             \, \dd W(s) \Big | ^{2q} \Big]     \notag  \\
 &  \leq &
          	C h^{2q-1}
          	 \int_t^{t+h} 
          	 \big \|f(X(s))- f(x) \big \|_{L^{2q}(\Omega, \R^d) }
          	 ^{2q}
          	 \, \dd s 
           +
              C h^{2q-1}
              \int_t^{t+h}
              \big \|f(x) - f(Y(t,x;t+h))
              \big \|_{L^{2q}(\Omega, \R^d) }
              ^{2q}
           	 \, \dd s      \notag \\
   & & \quad   +
             C h^{q-1}
          	 \int_t^{t+h}
          	   \big \|g(X(s)) - g(x) \big\|
          	   _{L^{2q}(\Omega, \R^{d \times m}) }
          	   ^{2q}
          	\, \dd s    \notag \\  
 & \leq  &
        C (1+|x|^{ 2q (4r+1) } )h^{2q},
\end{eqnarray}
as required. This thus finishes
the proof of the lemma.
\qed 

Based on Taylor's expansion and Lemma \ref{eq:norm-esti-one-step-BEM}, we are able to obtain the following local errors of the BEM \eqref{eq:numerical-approximation-BEM}.
\begin{lem}
\label{FBDE_one_step_estimation}
Under the same conditions of Lemma \ref{eq:norm-esti-one-step-BEM}, one can obtain the local errors of the BEM \eqref{eq:numerical-appro-one-step-BEM}
\begin{align}
\label{sheme:BEM_general_one_step_approximation}
     \bigg| \E \Big [ 
          \prod_{j=1}^s (\delta_{X,x})^{i_j}
       - 
           \prod_{j=1}^s (\delta_{Y,x})^{i_j} 
      \Big] \bigg|
   &  \leq 
        C (1+|x|^{8r+3})h^{2}, \,\, s=1,2,3,      \\
\label{sheme:norm-one-step-BEM}
      	\bigg\| \prod_{j=1}^{4}  
      	( \delta_{Y, x} )^{i_j} 
      	\bigg\|_{L^2 (\Omega; \R)}
   & \leq  
       C (1+|x|^{8r+4}) h^{2} .
\end{align}
\end{lem}
\textit{Proof of Lemma \ref{FBDE_one_step_estimation}.}
For $ s=1,$ we first note that 
\begin{eqnarray}
   & &   \big| \E \big[ 
            (\delta_{X,x})^{i_1}
         -  
             (\delta_{Y,x})^{i_1} 
       \big] \big|    
   =  
      \big| \E \big[
             (\delta_{X,x})^{i_1} 
         - 
             ( \delta_{Y_E,x})^{i_1} 
         -
             R_f^{i_1}
        \big]         \notag  \\
  & \leq &
      \big| \E \big[
             (\delta_{X,x})^{i_1} 
         - 
             \delta_{Y_E,x})^{i_1}
        \big] \big|
         +
     \bigg| \E \bigg[
             \int_t^{t+h}f^{i_1}
       \big( Y(t,x;t+h) \big)
         -
             f^{i_1}(x)
         \, \dd s \bigg] \bigg|.
\end{eqnarray}
The estimate for the first term has been done in
\eqref{eq:s1-error-one-step-Euler-method}, 
namely,
\begin{equation}
  \label{eq:s1-err-EM}
   \big| \E \big[
        (\delta_{X,x})^{i_1}
      - 
        ( \delta_{Y_E,x})^{i_1} 
    \big] \big|
      \leq 
        C (1+|x|^{4r+1})h^2.
\end{equation}
Next we bound the second term. 
Employing Taylor's expansion in the form of
\begin{eqnarray}
      f^{i_1} \big(Y(t,x;t+h) \big)   
& = &
       f^{i_1}(x) + 
       \sum^{d}_{j_1=1} 
        \frac {\partial f^{i_1}(x)}{\partial x_{j_{1}}}
       \big( 
         Y^{j_{1}}(t,x;t+h) 
      - 
         x^{j_{1}} \big)      \notag \\
 && \quad  + 
         \sum^{d}_{j_{1},j_{2}=1} 
         \int_0^1
       \frac{\partial ^{2} f^{i_1}
      \big( x + \theta \big( Y(t,x;t+h)-x 
       \big)\big)}
       {\partial x_{j_{1}} \partial x_{j_{2}}}
       (1-\theta) 
       \, \dd \theta      \notag \\
 & & \quad  
   \cdot  \big(
           Y^{j_1}(t,x;t+h) -  x^{j_1}  \big) 
        \big( Y^{j_2}(t,x;t+h) - x^{j_2} \big),
\end{eqnarray}
and using the assumption on the function $ f $ in \eqref{ass:drift-f} and the H\"older inequality, one can infer 
\begin{eqnarray}
  & &
   \Big|  \E \Big[
      \int_t^{t+h}f^{i_1}
       \big(Y(t,x;t+h) \big)
        - f^{i_1}(x)
          \, \dd s
           \Big] \Big|     \notag  \\
 & \leq &
       \sum^{d}_{j_1=1}
         \Big|  \E \Big[
          \int_t^{t+h} 
           \frac {\partial f^{i_1}(x)}
                 {\partial x_{j_{1}}}
             \big(Y^{j_1}(t,x;t+h)-x^{j_1}
              \big) 
               \, \dd s \Big] \Big|   \notag  \\
 & & \quad
        + 
           \sum^{d}_{j_{1},j_{2}=1}
           \Big| \E \Big[
            \int_t^{t+h} \int_0^1
            \frac{\partial ^{2} f^{i_1}(x + \theta (Y(t,x;t+h)-x))}
                 {\partial x_{j_{1}} \partial x_{j_{2}}} (1-\theta) 
             \, \dd \theta   
              \,\dd s        \notag \\
 & & \quad \cdot
       \big( Y^{j_{1}}(t,x;t+h) - x^{j_{1}}\big )
       \big( Y^{j_{2}}(t,x;t+h) - x^{j_{2}} \big )
        \Big] \Big|      \notag \\
 & \leq &
      C (1+|x|^{2r})
        \Big| \E \Big[ 
         \int_t^{t+h}
          \big(Y^{j_1}(t,x;t+h)-x^{j_1}
           \big)
           \, \dd s
            \Big] \Big|     \notag \\
 & & \quad
       + 
       C (1+|x|^{2r-1}) h \,
        	\| (\delta_{Y,x})^{j_1}
        	\|_{L^{2} (\Omega; \R) }
          	\| (\delta_{Y,x})^{j_2}
          	\|_{L^{2} (\Omega; \R) }    \notag  \\
 & \leq &
          C (1+|x|^{6r+1})h^2.
\end{eqnarray}
This together with \eqref{eq:s1-err-EM}
implies
 \begin{equation}
   \label{eq:BEM_one_step_error_s1}
   \big| \E \big[
     (\delta_{X,x})^{i_1}
  - 
     (\delta_{Y,x})^{i_1}
   \big] \big|
    \leq 
     C (1+|x|^{6r+1})h^2.
   \end{equation}
For $ s=2 $, we first note that
\begin{eqnarray}
 & &
     \big| \E \big[ 
        (\delta_{X,x})^{i_1}(\delta_{X,x})^{i_2}
      - 
        (\delta_{Y,x})^{i_1}(\delta_{Y,x})^{i_2}
       \big] \big|     \notag \\
 & \leq &
     \big|\E \big[
        (\delta_{X,x})^{i_1} ((\delta_{X,x})^{i_2}
      - 
        (\delta_{Y,x})^{i_2}) \big]\big|
      + 
        \big| \E \big[
        ((\delta_{X,x})^{i_1}
      - 
        (\delta_{Y,x})^{i_1})
        (\delta_{Y,x})^{i_2}
        \big]\big|         \notag  \\
 & = &
      \big| \E \big[
        (\delta_{X,x})^{i_1}
        ((\delta_{X,x})^{i_2}
      -
        ( \delta_{Y_E,x})^{i_2} - R_f^{i_2})
      \big]\big|
      + 
        \big| \E \big[
        ((\delta_{X,x})^{i_1}
        - ( \delta_{Y_E,x})^{i_1}-R_f^{i_1})
        (\delta_{Y,x})^{i_2}
      \big]\big|       \notag \\
 & \leq &
      \big|\E \big[
        (\delta_{X,x})^{i_1}((\delta_{X,x})^{i_2}
      - 
       ( \delta_{Y_E,x})^{i_2})
      \big]\big|
     +
      \big| \E \big[
      ((\delta_{X,x})^{i_1}
      -
    (\delta_{Y_E,x})^{i_1})(\delta_{Y,x})^{i_2}
      \big] \big|     \notag \\
& & \quad +
      \big|\E \big[
            (\delta_{X,x})^{i_1}R_f^{i_2}
              \big]\big|
     +
      \big|\E \big[ 
           R_f^{i_1}(\delta_{Y,x})^{i_2}
              \big]\big|.
\end{eqnarray}
By \eqref{eq:esti_I_2},
it follows that
\begin{equation}
  \label{eq:err-exact-EM}
    \big|\E \big[
        (\delta_{X,x})^{i_1}
        ((\delta_{X,x})^{i_2}
      - 
      ( \delta_{Y_E,x})^{i_2})
       \big]\big|
  \leq
     C (1+|x|^{6r+2})h^{2}.
\end{equation}
The estimate for the second is obtained in the same way. 
Utilizing the H\"older inequality,
\eqref{eq:norm-estimate-one-step-exact-solution}
and
\eqref{eq:norm-one-step-Rf}
yields
\begin{equation}
  \label{eq:esti-norm-exact-R}
     \big|\E \big[
        (\delta_{X,x})^{i_1} R_f^{i_2}
       \big]\big| \leq
       	\| (\delta_{X,x})^{i_1} \|_{L^{2} (\Omega; \R) } \,
       \| R^{i_2}_f \|_{L^{2} (\Omega; \R) }
    \leq
   C(1+|x|^{6r+2})h^{2}.
\end{equation}
Following a similar way as \eqref{eq:esti-norm-exact-R},
we acquire
\begin{equation}
 \label{eq:esti-norm-BEM-R}
   \big|\E \big[
          R_f^{i_1} (\delta_{Y,x})^{i_2}
     \big] \big|
 \leq
    C(1+|x|^{6r+2})h^{2}.
\end{equation}
Collecting
\eqref{eq:err-exact-EM},
\eqref{eq:esti-norm-exact-R},
and
\eqref{eq:esti-norm-BEM-R}
enables us to arrive at
\begin{equation}
 \label{eq:BEM_one_step_error_s2}
	\big| \E \big[ 
	  (\delta_{X,x})^{i_1} (\delta_{X,x})^{i_2}
	- 
	  (\delta_{Y,x})^{i_1} (\delta_{Y,x})^{i_2}
	 \big]\big|
	  \leq 
	   C(1+|x|^{6r+2})h^{2}.
\end{equation}
Finally, for $ s=3,$ 
we make a decomposition as follows:
\begin{eqnarray}
& &
  \big|\E \big[ 
      (\delta_{X,x})^{i_1} (\delta_{X,x})^{i_2} (\delta_{X,x})^{i_3}
     - 
      (\delta_{Y,x})^{i_1} (\delta_{Y,x})^{i_2} (\delta_{Y,x})^{i_3}\big]\big|
      \notag  \\
& \leq &
        \big| \E\big[
         (\delta_{X,x})^{i_1}((\delta_{X,x})^{i_2} (\delta_{X,x})^{i_3}
       - 
         (\delta_{Y,x})^{i_2}(\delta_{Y,x})^{i_3}) \big]\big|
      +
        \big| \E\big[
        ((\delta_{X,x})^{i_1}-(\delta_{Y,x})^{i_1})
           (\delta_{Y,x})^{i_2} (\delta_{Y,x})^{i_3}\big]\big|
           \notag  \\
 & \leq &
        \big|\E \big[
             (\delta_{X,x})^{i_1}(\delta_{X,x})^{i_2}
             ((\delta_{X,x})^{i_3}-(\delta_{Y,x})^{i_3})
          \big]\big|
      +
        \big|\E \big[
                 (\delta_{X,x})^{i_1}(\delta_{Y,x})^{i_3}
                 ((\delta_{X,x})^{i_2}-(\delta_{Y,x})^{i_2})
          \big]\big|      \notag \\
 & &  \quad  +
        \big|\E\big[
             ((\delta_{X,x})^{i_1}-(\delta_{Y,x})^{i_1})(\delta_{Y,x})^{i_2}
            (\delta_{Y,x})^{i_3}  \big]\big|.
\end{eqnarray}
Using
\eqref{eq:norm-estimate-one-step-exact-solution},
\eqref{eq:norm-one-step-BEM}, 
\eqref{eq:norm-error-one-step-solution-BEM}
and the H\"older inequality yields
\begin{equation}
  \label{eq:decom_1_BEM_one_step_error}
  \big| \E \big[
           (\delta_{X,x})^{i_1}(\delta_{X,x})^{i_2}
           ((\delta_{X,x})^{i_3}-(\delta_{Y,x})^{i_3})
           \big]\big|
  \leq C (1 + |x|^{8r+3}) h^2,
\end{equation}
\begin{equation}
   \label{eq:decom_2_BEM_one_step_error}
  \big|\E \big[(\delta_{X,x})^{i_1}(\delta_{Y,x})^{i_3}
  ((\delta_{X,x})^{i_2}-(\delta_{Y,x})^{i_2})\big]\big|
  \leq C (1 + |x|^{8r+3}) h^2,
\end{equation}
and
\begin{equation}
  \label{eq:decom_3_BEM_one_step_error}
  \big|\E \big[ ((\delta_{X,x})^{i_1}-(\delta_{Y,x})^{i_1})
  (\delta_{Y,x})^{i_2}(\delta_{Y,x})^{i_3}\big]\big|
  \leq C (1 + |x|^{8r+3}) h^2.
\end{equation}
Putting \eqref{eq:BEM_one_step_error_s1},
\eqref{eq:BEM_one_step_error_s2} and
\eqref{eq:decom_1_BEM_one_step_error}-
\eqref{eq:decom_3_BEM_one_step_error} together, we obtain
\begin{equation}
  \Big| \E \Big[ \prod_{j=1}^s (\delta_{X,x})^{i_j} \Big] 
  -
   \E \Big[ \prod_{j=1}^s (\delta_{Y,x})^{i_j} \Big] \Big|
   \leq  C (1+|x|^{8r+3}) h^2,
   \,\,
   s = 1,2,3.
\end{equation}
Similar to the proof of \eqref{sheme:solution_one-step approximation2}, one can employ the H\"older inequality and 
\eqref{eq:norm-one-step-BEM},
to deduce
\begin{equation}
	\bigg\| \prod_{j=1}^{4}  
	 ( \delta_{Y, x} )^{i_j} 
	 \bigg\|_{L^2 (\Omega; \R)}
  \leq   
	  \prod_{j=1}^{4}
  	   \Big\|   
	    ( \delta_{Y, x} )^{i_j} 
	    \Big\|_{L^8 (\Omega; \R)}
    \leq
       C (1+|x|^{8r+4}) h^2,
 \end{equation}
 which finishes the proof of Lemma 
\ref{FBDE_one_step_estimation}.
\qed

Based on the local error and Corollary
\ref{cor-fundamental-weak-convergence}, we are able to obtain the global
weak convergence rate of the BEM
\eqref{eq:numerical-approximation-BEM}.
\begin{thm}
 \label{thm:BEM-global-conver-rate}
Suppose Assumption 
\ref{ass:coefficient_function_assumption} hold with $p=1$.
Then the backward Euler method  \eqref{eq:numerical-approximation-BEM} has a global weak convergence rate of order $ 1 $, namely, for $ \varphi \in \mathbb{H}^{4} $
  \begin{equation}
   \big |\E \big [
	\varphi \big (
	X(t_0,X_0;T)
	\big ) \big ]
	-\E \big[
	\varphi \big (
	Y( t_0,Y_0;T )
	\big )\big ] \big |
	\leq 
	C 
	\big( 1+ |X_0|^{ (8r+4+\chi) } \big) h.
	\end{equation} 
\end{thm}
\section{Numerical experiments}
\label{sec:numer-results}
In this section, we present some numerical experiments to illustrate the previous theoretical findings.
In the following simulations, we used the Mersenne Twister random generator with seed $100$. 
\textcolor{black}{
Newton's method is used to solve the nonlinear algebraic equations at each step of the implicit schemes with the two successive iterative solutions smaller than $10^{-6}$ in the $l^2$-norm.}
\par
We test several numerical schemes for two model problems. The first one satisfies Assumption \ref{ass:coefficient_function_assumption} (non-globally Lipschitz both drift and diffusion). The second example has non-globally Lipschitz drift and global Lipschitz diffusion. The aim of the tests is to compare the weak convergence rates of the numerical methods proposed in subsection \ref{ssec:example-mod-euler}.
\par
Along this section, we consider four different test functions:
 \begin{equation}
 \label{eq:test_function}
     \varphi(x) 
  = x, \, \, x^2, \,\, \cos(x), \,\, e^{-x^2}.
 \end{equation}
The reference values of the models are computed based on a Monte Carlo method combined with the scheme \eqref{eq:numerical-approximation-BEM} and $ h_{ref} =2^{-14} $:
 \begin{equation}
     \E \big[ \varphi(X_T)\big]
     \thickapprox
    \frac{1}{M} \sum_{i=1}^{M}
    \varphi\big( X_T(\omega_i,h_{ref} )\big).
 \end{equation}
\begin{example}
\label{ex:numerical_ex_nonlinear_diff}
As the first test model, we consider the following model
\begin{equation}
  \label{eq:non_linear_model}
	\left\{
	\begin{aligned}
		\dd X(t) & = \big(1- X^5(t) + X^3(t)\big) \, \dd t 
+ \big( \tfrac{1}{10} X^2(t) + 2 \big) \, \dd W(t),
		\quad t \in (0,T],
		\\
		X(0) & = x.
	\end{aligned}\right.
\end{equation}
\end{example}
{\color{black}
The conditions in Assumption \ref{ass:coefficient_function_assumption} are fulfilled
with $r=2$ and $\rho=2$.

In this example,  we test the weak errors of the backward Euler method (BEM) \eqref{eq:numerical-approximation-BEM}, the fully tamed Euler (FTE1) \eqref{eq:numerical-approximation-FTE1}, the fully tamed Euler (FTE2) \eqref{eq:numerical-approximation-FTE2}, the modified Euler scheme (MES) \eqref{eq:numerical-approximation-MES}, the balanced scheme (BS) \eqref{eq:numerical-approximation-BS} and the balanced-type scheme (BTS) \eqref{eq:numerical-approximation-BTS}.}
\par
For the first presented numerical experiments, we consider a unit terminal time $ T=1 $, the initial condition $ x = 2, 8 $, and the time step sizes $ h = 2^{-k}, k= 6,7,8,9,10.$ In addition, the quantity $\E \big[ \varphi(Y_T)\big]$ is estimated by a Monte Carlo approximation, involving $ M = 3 \times 10^6 $ independent trajectories.
The statistical errors for the Monte Carlo approximations are computed with $95\%$ confidence \cite[subsection 2.5.1]{Zhang2017numerical} and the resulting statistical errors are 
10 times smaller than the reported weak errors.
\par
 
\par
\begin{figure}
    \includegraphics[width=0.5\linewidth, height=0.3\textheight]{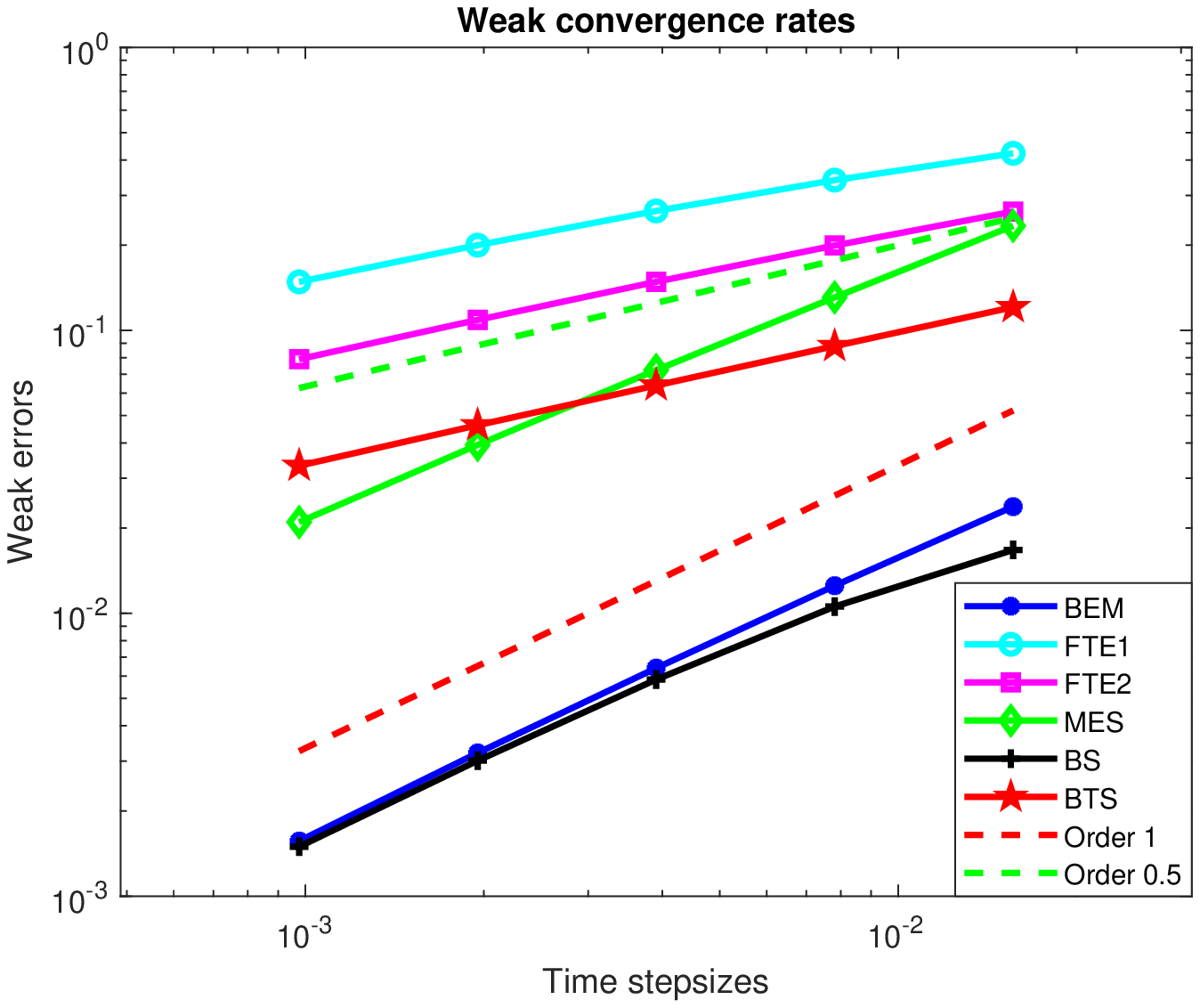}
    \includegraphics[width=0.5\linewidth, height=0.3\textheight]{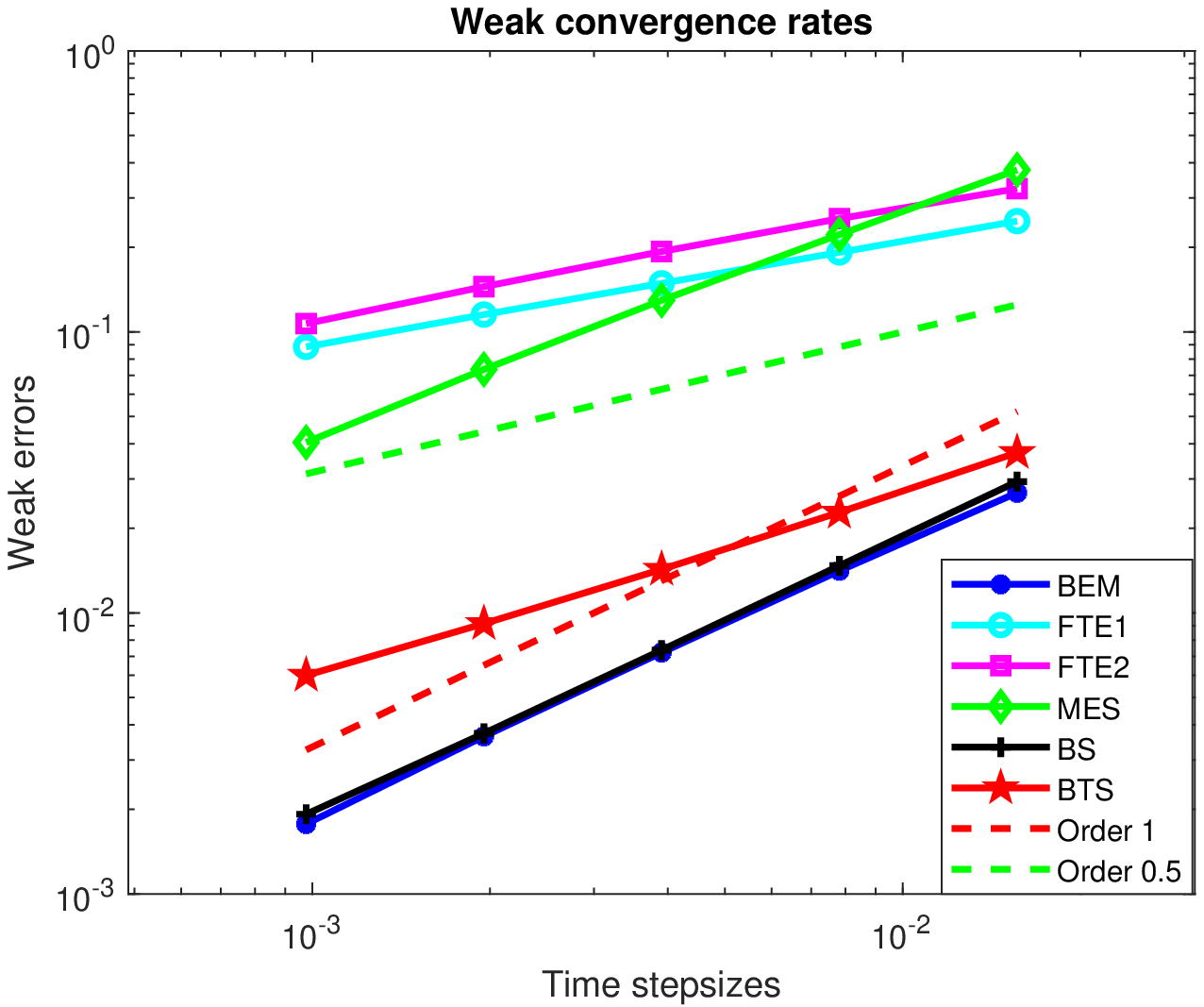}
	\caption
	    {Weak convergence rates for the stochastic differential equation	  \eqref{eq:non_linear_model}
	    	for $ \varphi(x) = x $ (Left) 
	    	and $ \varphi(x) = x^2 $ (Right).
	     }
\label{fig:nonlinear_diff_x_x_2!}	    
\end{figure}
\begin{figure}
	\includegraphics[width=0.5\linewidth, height=0.3\textheight]{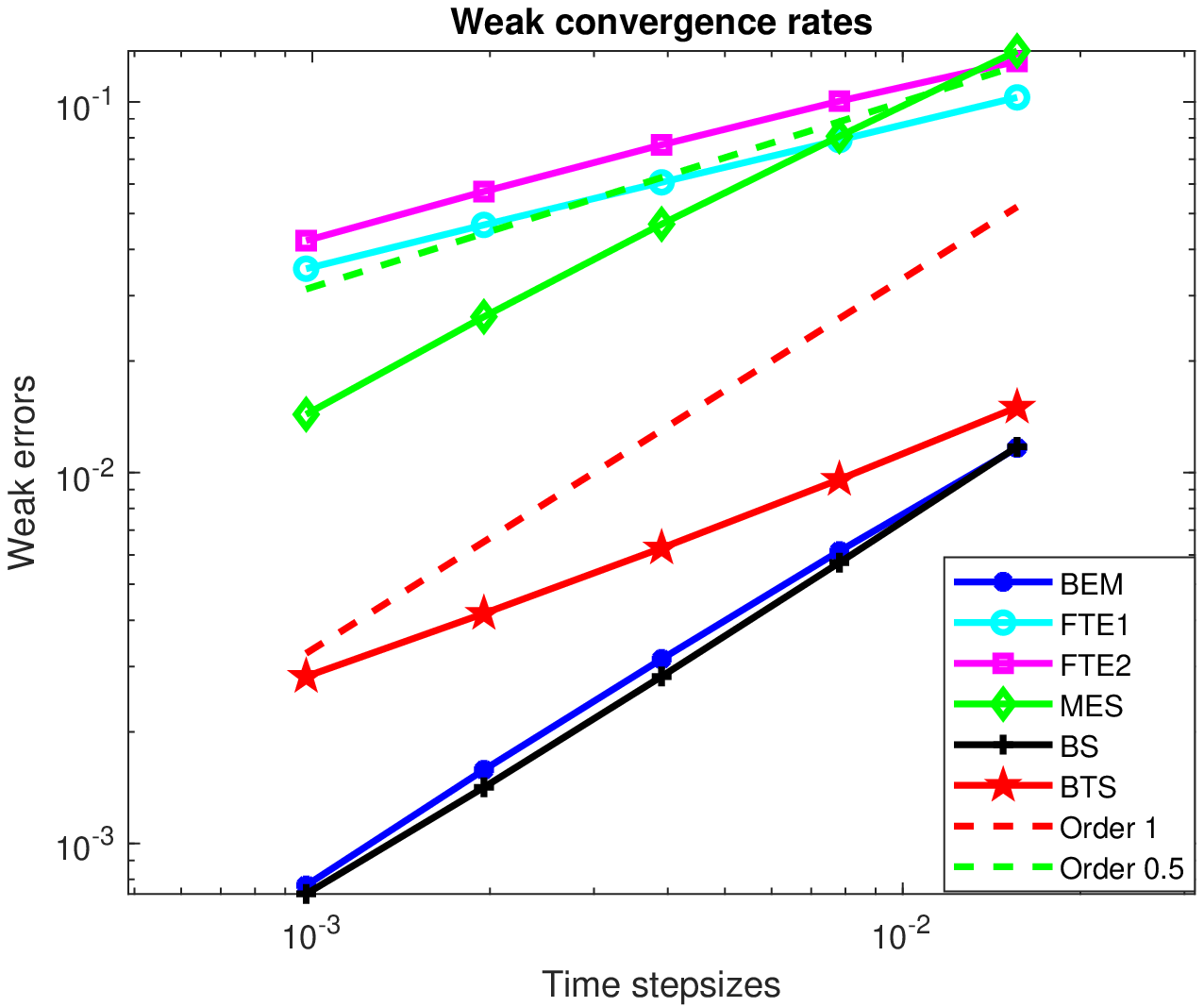}
	\includegraphics[width=0.5\linewidth, height=0.3\textheight]{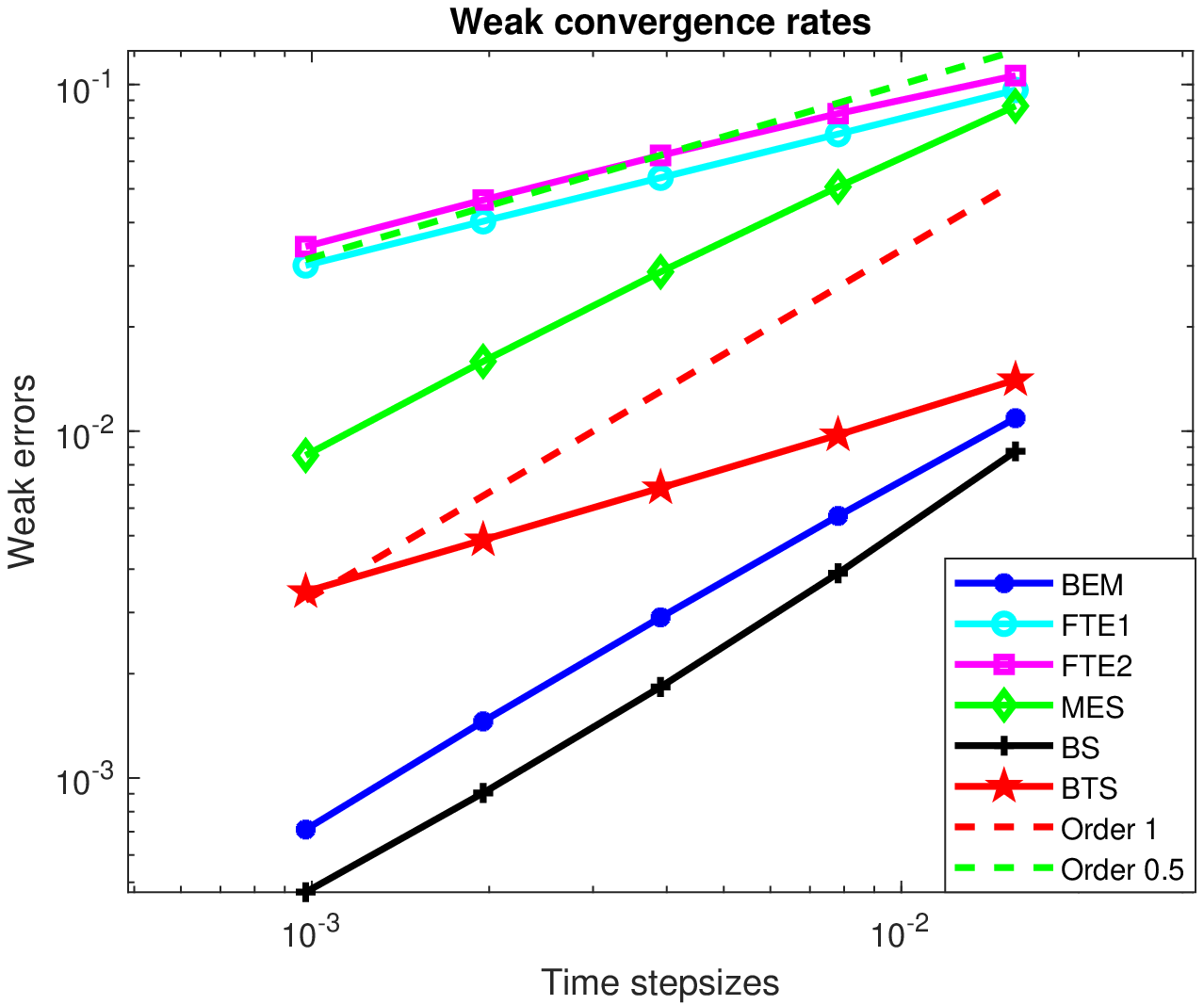}
	\caption
	    {Weak convergence rates for the the stochastic differential equation
	   	\eqref{eq:non_linear_model}
	    	for $ \varphi(x) = \cos(x) $ (Left) 
	    	and $ \varphi(x) = e^{-x^2} $ (Right).
	     }
  \label{fig:nonlinear_diff_cosx_exp_x_2!}
\end{figure}

In Figure \ref{fig:nonlinear_diff_x_x_2!} and Figure \ref{fig:nonlinear_diff_cosx_exp_x_2!}, we take $x=2$ and  present the weak errors of six different numerical methods with the test functions \eqref{eq:test_function} against five different stepsizes on a log-log scale.
As expected, the weak convergence rates of the MES \eqref{eq:numerical-approximation-MES}, the BEM \eqref{eq:numerical-approximation-BEM}, and the BS \eqref{eq:numerical-approximation-BS} are close to $1$, whereas the FTE1 \eqref{eq:numerical-approximation-FTE1}, the FTE2 \eqref{eq:numerical-approximation-FTE2}, and the BTS \eqref{eq:numerical-approximation-BTS} give convergence rates close to $\frac{1}{2}$.
{\color{black} When $x=8$, our numerical experiments show that the tamed-type methods including the FTE1 \eqref{eq:numerical-approximation-FTE1}, the FTE2 \eqref{eq:numerical-approximation-FTE2}, and the MES \eqref{eq:numerical-approximation-MES} can not produce satisfactory accuracy with the considered time step sizes $ h = 2^{-k}, k= 6,7,8,9,10$, while the BS \eqref{eq:numerical-approximation-BS}, the BTS \eqref{eq:numerical-approximation-BTS} and the BEM \eqref{eq:numerical-approximation-BEM} always give high accuracy approximations. In this case, we also implement the Euler-Maruyama (EM) method for the test problem and find that  all ($100\%$) the EM trajectories explode, even for very small step-size $h = 2^{-10} \approx 0.001 $.
}
\vspace{0.2cm}
\begin{example}
\label{ex:numerical_ex_linear_diff}
Consider the stochastic
FitzHugh-Nagumo (FHN) model \cite{buckwar2022splitting} in the form of:
\begin{equation}
	\label{eq:stochastic_FHN_model}
	\begin{split}
		&  \left(\begin{array}{c}
			\dd X_{1}(t) \\ \dd X_{2}(t) \end{array}\right)
		  = \left(\! \begin{array}{c}
		 	X_{1}(t)-X_{1}^3(t) - X_{2}(t)
		 	\\  X_{1}(t) - X_{2}(t) + 1 
	 	\end{array} \! \right) \dd t
 	+  \left(\! \begin{array}{cc}
 	      X_1(t) + 1 & 
 		0 \\ 0 & 
 		 X_2(t) + 1
 		\end{array} \! \right)
 	   \dd  W(t)
	\end{split}
\end{equation}
for $ t \in (0,T]$ and $ X(0)  = [0,0] $, with solution $ X(t) := \big( X_1(t) , X_2(t) \big)^{T} $
for $ t \in [0 , T] $,
where $ W(t) : = \big( W_1(t), W_2(t) \big)^{T} $
is a two dimensional Brownian motion.
\end{example}
\par
In this example, conditions in Assumption \ref{ass:coefficient_function_assumption} are fulfilled
with $r=1$ and $\rho=1$.
Here we test weak convergent rates of the BEM \eqref{eq:numerical-approximation-BEM}, the FTE1 \eqref{eq:numerical-approximation-FTE1}, the FTE2 \eqref{eq:numerical-approximation-FTE2}, the MES \eqref{eq:numerical-approximation-MES}, the DTE scheme \eqref{eq:numerical-approximation-DTE}, and the BS \eqref{eq:numerical-approximation-BS}.
\par

For the second  numerical experiment, we consider a unit terminal time $ T=1 $ and the time step $ h = 2^{-k}, k= 7,8,9,10,11$, and $M = 10^6 $ independent trajectories in the Monte Carlo approximation.
The statistical errors for the Monte Carlo approximations turn out to be negligible as the Monte Carlo error computed with 95\% confidence is at least 10 times smaller than the weak errors.
\par
\begin{figure}	\includegraphics[width=0.5\linewidth, height=0.3\textheight]{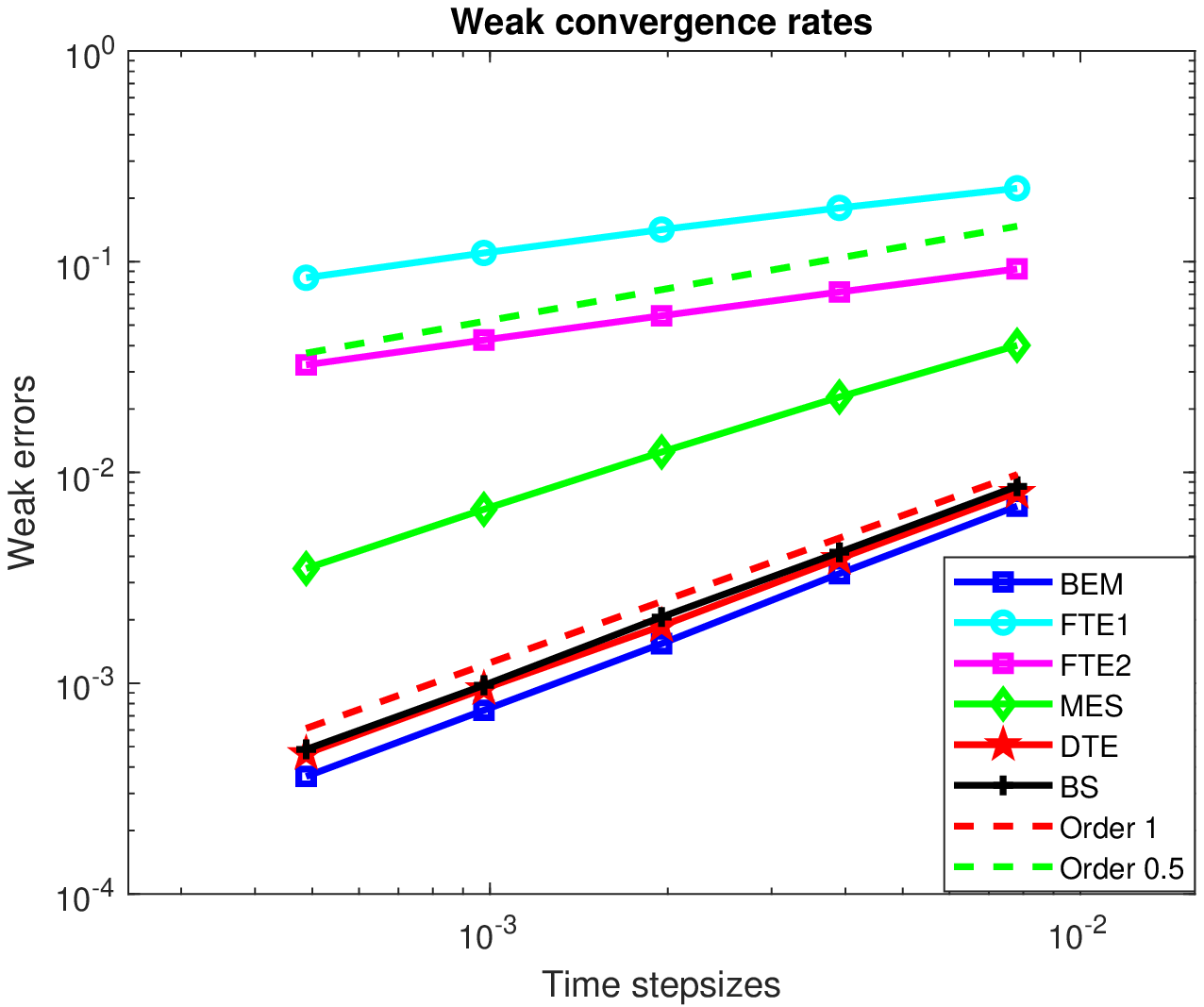}
	\includegraphics[width=0.5\linewidth, height=0.3\textheight]{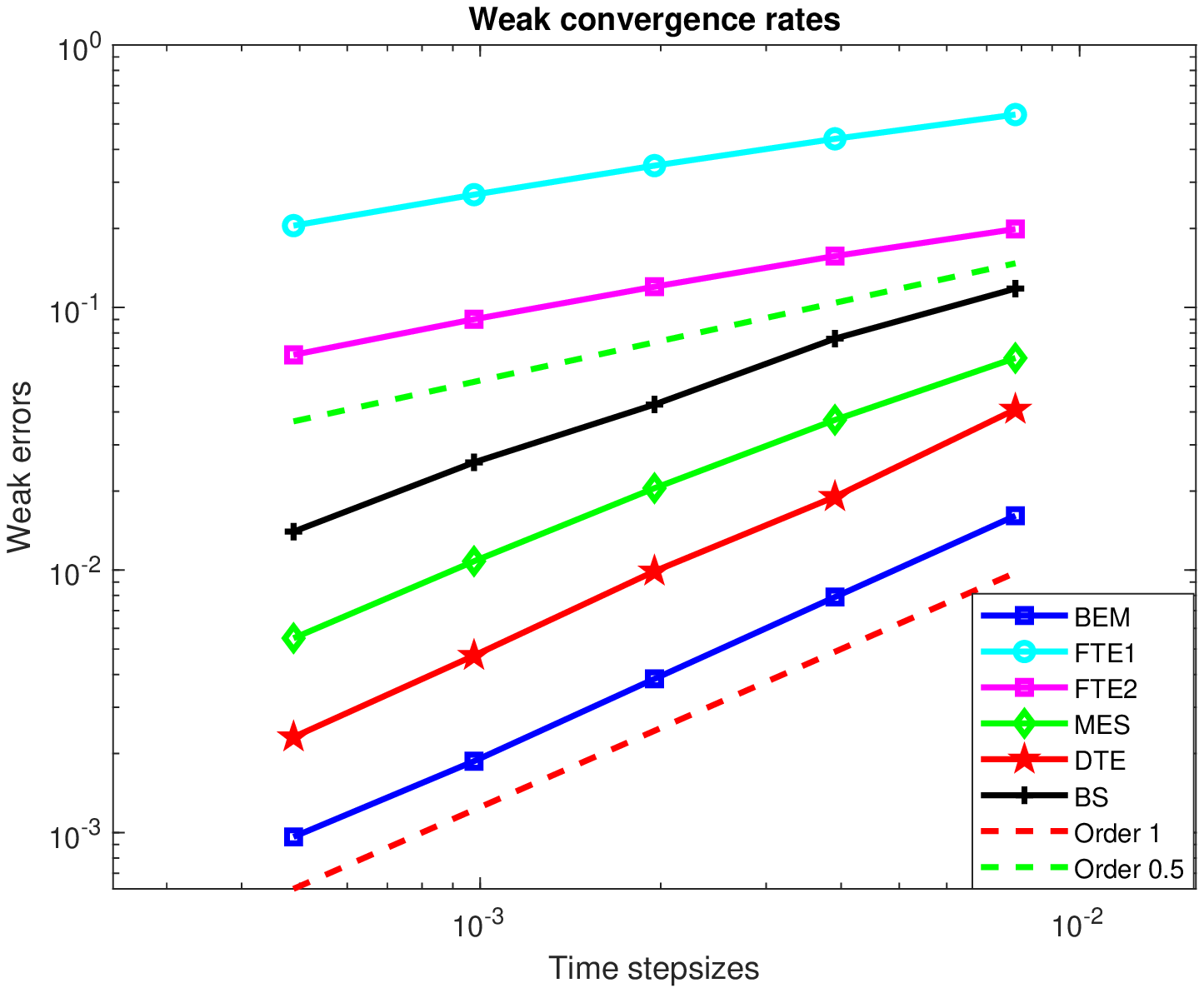}
	\caption{Weak rates for stochastic FHN model with
		  $ \varphi(x) = x$ (Left)
		  and
		  $ \varphi(x) = x^2 $ (Right).}
 \label{fig:linear_diff_x_x_2!}
\end{figure}
\begin{figure}[h]	\includegraphics[width=0.5\linewidth, height=0.3\textheight]{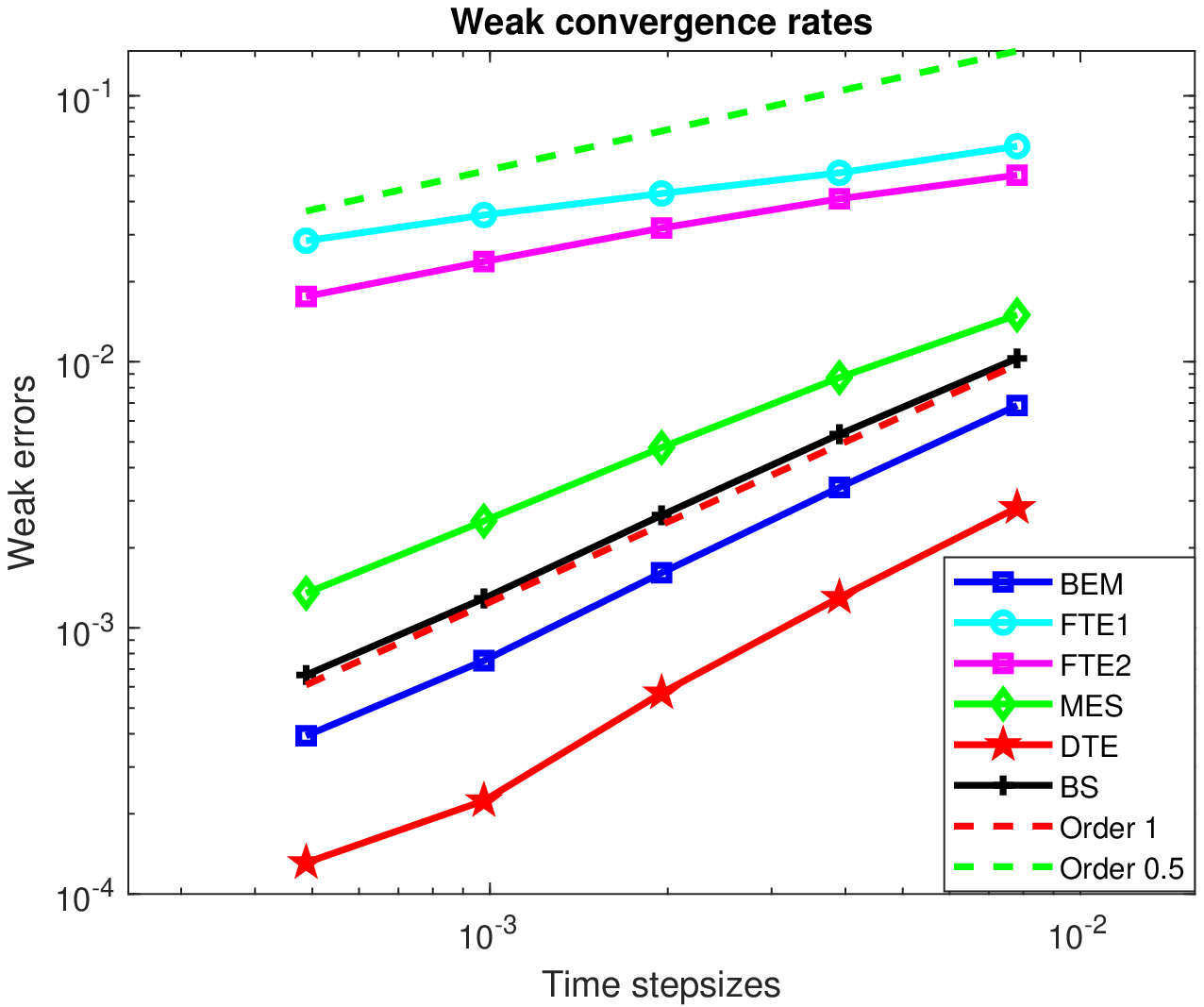}
\includegraphics[width=0.5\linewidth, height=0.3\textheight]{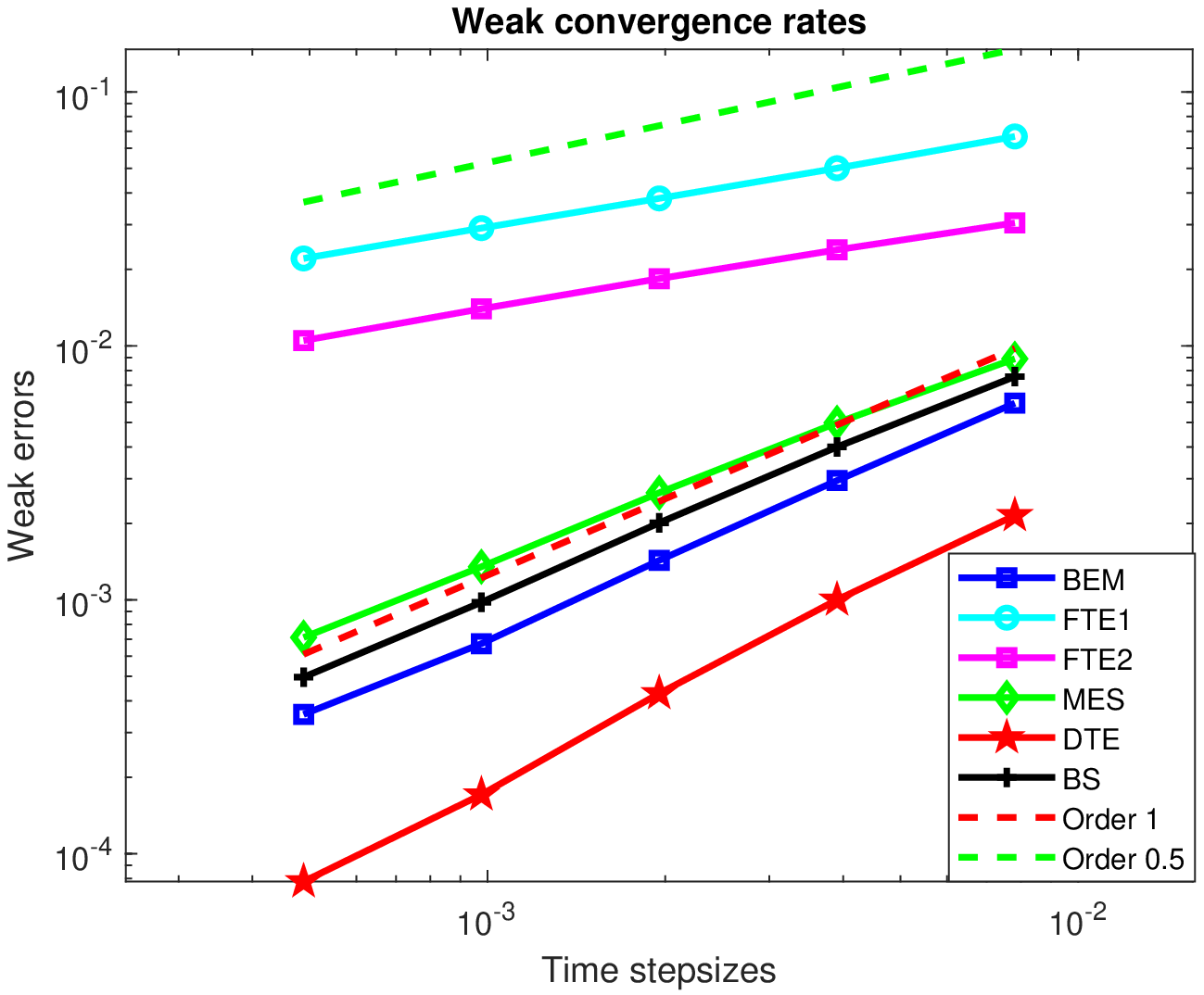}
	\caption{Weak rates for stochastic FHN model with
		  $ \varphi(x) = \cos(x) $ (Left)
		  and
		  $ \varphi(x) = e^{-x^2} $ (Right).}
 \label{fig:linear_diff_cosx_exp_x_2!}
\end{figure}

From Figure \ref{fig:linear_diff_x_x_2!} and Figure \ref{fig:linear_diff_cosx_exp_x_2!},
we observe that weak convergence
rates of 
\eqref{eq:numerical-approximation-MES}, \eqref{eq:numerical-approximation-DTE}, \eqref{eq:numerical-approximation-BS},
\eqref{eq:numerical-approximation-BEM}
are $ 1 $.
For the schemes \eqref{eq:numerical-approximation-FTE1}
and
\eqref{eq:numerical-approximation-FTE2}, we observe a 
weak convergence of order 
$ \frac{1}{2} $.

\vskip6mm
\bibliographystyle{plain}
\bibliography{main}


\appendix

\section{Proof of Theorem \ref{thm:fundamental_weak_convergence}}
\label{append:proof-fund}
%
The proof exploits the idea of the proof of \cite[Theorem 2]{Milstein1985Weak} and \cite[Theorem 2.1]{Milstein2004stochastic}
in the global Lipschitz case. Recall $Y_0 = X_0$ and
\begin{equation}
\label{eq:X-flow-property}
X \big (t , x ; T \big )
=X \big (s, X_{t, x }(s) ; T \big ),
\quad
0 \leq t \leq s \leq T.
\end{equation}
Noting $t_N = T$ and recalling the definition \eqref{eq:u-defn} yield
\begin{equation}
\E \big [
  	 \varphi \big (
    	X(t_0,X_0;T)
  	    \big ) \big ]
  	     -\E \big[
  	      \varphi \big (
     	Y(t_0,Y_0; t_N)
  	   \big )\big ]
=
\E
[ u ( t_0, X_0) ]
-
\E 
[
u (t_N, Y_N )
].
\end{equation}
Again, by the definition \eqref{eq:u-defn} and \eqref{eq:X-flow-property},
\begin{eqnarray}
u ( t, x ) & =&
\E \big[ \varphi (X(t, x; T) ) \big]
=
\E \big[ \varphi (X( s, X_{t, x} ( s ); T) ) \big] \notag
\\
& =&
\E \big[ 
\E \big( \varphi (X( s, X_{t, x} ( s ); T) ) | \sigma ( X_{t, x} ( s ) ) \big) \big] 
 =
\E [ u ( s, X_{t, x} ( s ) ) ],
\quad
\text{ for }
0 \leq t \leq s \leq T. 
\end{eqnarray}
By taking 
  $ t = t_i, $
  $ s = t_{i+1}, $
a conditional version of such an equality 
  $ u ( t_i, Y_i )
             =
               \E \big[ 
                       u \big( t_{i+1}, X_{t_i, Y_i} ( t_{i+1} 
                         \big)
                   | \sigma ( Y_i )  
                  \big],
  $
implies             
\begin{equation}
\label{eq:u-expectation}
\E \big[ u(t_i , Y_i)  \big]
=
\E \big[ u(t_{i+1} , X (t_i, Y_i; t_{i+1}) )  \big] ,
\quad
i = 0, 1,..., N-1.
\end{equation}
Using this repeatedly, we obtain
\begin{eqnarray}
 \label{eq:weak-error-sum}
		\E \big[ u(t_0,X_0) \big]
		 & = &
		\E \big[ u(t_1 , X( t_{1} ) ) \big]
		- \E \big[ u(t_1 , Y_1)  \big]
		+ \E \big[ u( t_{2}, X(t_1,Y_1; t_2) )  \big] \notag \\
		&= &
		\E \big[ u(t_1 , X( t_0, X_0; t_{1} ) ) \big]
		- \E \big[ u(t_1 , Y_1)  \big] 
		\notag \\
		&& \quad +
		\E \big[ u(t_2 , X( t_{1}, Y_1; t_2 ) ) \big]
		- \E \big[ u(t_2 , Y_2)  \big]
		+ \E \big[ u( t_{3}, X(t_2,Y_2; t_3) )  \big] \notag \\
	   &= & \cdots  \notag \\
	   &= &
		\sum_{i=0}^{N-2}
		\Big( \E \big[ u(t_{i+1} , X( t_{i}, Y_{i}; t_{i+1} ) ) \big]
		- \E \big[ u(t_{i+1} , Y_{i+1})  \big] \Big) 
		 +
		\E \big[ u( t_{N-1},Y_{N-1} ) \big].		
\end{eqnarray}
By
\eqref{eq:u-expectation}, we get
\begin{equation}
\E \big[ u( t_{N-1},Y_{N-1} ) \big]
		-
		\E \big[ u(t_N,Y_N) \big] 
=
\E \big[ u(t_{N} , X (t_{N-1}, Y_{N-1}; t_N) )  \big]
-
\E \big[ u(t_N,Y_N) \big],
\end{equation}
and thus one can derive from \eqref{eq:weak-error-sum} that
\begin{eqnarray}
	\label{eq:error_test_function}
	&&
		\E \big[ u(t_0,X_0) \big] - \E \big[ u(t_N,Y_N) \big]  \notag \\
		& = &
		\sum_{i=0}^{N-1}
		\Big( \E \big[ u(t_{i+1} , X( t_{i}, Y_{i}; t_{i+1} ) ) \big]
		- \E \big[ u(t_{i+1} , Y_{i+1})  \big] \Big)
		\notag \\
		& = & 
		\sum_{i=0}^{N-1}
		\Big( \E \big[ u(t_{i+1} , X( t_{i}, Y_{i}; t_{i+1} ) ) \big]
		- \E \big[ u(t_{i+1} , Y( t_{i}, Y_{i}; t_{i+1} ) )  
		\big] \Big).
\end{eqnarray}
Finally, combining a conditional version of \eqref{eq:u-local-error} with \eqref{eq:error_test_function} yields
\begin{eqnarray}
&&
\big | \E \big [
  	 \varphi \big (
    	X(t_0,X_0;T)
  	    \big ) \big ]
  	     -\E \big[
  	      \varphi \big (
     	Y(t_0,Y_0; t_N)
  	   \big )\big ] \big |
 \notag \\
& \leq &
\sum_{i=0}^{N-1}
		\Big| \E \big[ u(t_{i+1} , X( t_{i}, Y_{i}; t_{i+1} ) ) \big]
		- \E \big[ u(t_{i+1} , Y( t_{i}, Y_{i}; t_{i+1} ) )  \big] \Big|
 \notag \\
& = &
\sum_{i=0}^{N-1}
		\Big| \E \Big[  \E \big( u(t_{i+1} , X( t_{i}, Y_{i}; t_{i+1} ) ) - u(t_{i+1} , Y( t_{i}, Y_{i}; t_{i+1} ) ) | \sigma(Y_i) \big) \Big] \Big|
 \notag \\
& \leq &
\sum_{i=0}^{N-1}
\E [ | \nu ( Y_i ) | ] h^{p+1}
\leq
\Upsilon (X_0) T h^p,
\end{eqnarray}
which validates the desired assertion \eqref{eq:thm-convergence-rate}.
\qed
\section{Proof of Lemma \ref{lem:equivent_condition}}
\label{append:proof-one-step}
\textit{Proof}.
Since $\varphi (x)$ is differentiable and the solution $ X(t, x; T) $  of the SDE \eqref{eq:SDE-X(t,x;T)} is 
mean-square differentiable with respect to the initial data $x$ by assumption, 
$u (t, x) $ is differentiable with respect to $x$ and
\begin{equation}
\tfrac{\partial u ( t, x )}{ \partial x_{i} } 
= 
\tfrac{\partial
		\E [ \varphi (X(t,x;T) ) ] }{\partial x_{i}
		}
	=
	\E
	[
	D \varphi ( X ( t,x;T) )
	\tfrac{\partial
		X ( t,x;T) } {\partial x_{i}
		}
	]
	=
	\E 
	\Big[
	\sum_{j = 1}^d \partial_{x_j} \varphi  (X(t,x;T) ) 
	\tfrac{\partial
		X^j ( t,x;T) } {\partial x_{i}
		}
	\Big].
\end{equation}
Noting further that the partial derivative of $ \varphi $ belongs to $\mathbb{H}$ 
by assumption, 
one can use \eqref{eq:ass-Derivative-X}
and the H\"older inequality to deduce
\begin{equation} \label{eq:1st-partial-u}
| \tfrac{\partial u ( t, x )}{ \partial x_{i} } |
\leq
\| D \varphi (X(t,x;T) ) \|_{L^2 ( \Omega, \R^d ) } 
\| \tfrac{\partial
		X(t,x;T)} {\partial x_{i}
		} \|_{L^2 ( \Omega, \R^d ) }
\leq
C ( 1 + | x |^{\chi} ).
\end{equation}
For higher order derivatives, we recall that $ \varphi $'s partial derivatives of 
order up to $ 2p + 2 $ inclusively, belong to $\mathbb{H}$ and 
the solution of the SDE is $ 2p + 2 $ times mean-square differentiable with respect to the initial data $x$.
These facts together with the chain rule enable us to similarly show that 
$u (t, x) $ is $ 2p + 2 $ times differentiable with respect to $x$ and
\begin{equation}
\label{eq:kth-partial-u}
\big|
\tfrac{\partial ^{k} u ( t + h, x)}{\partial x_{i_{1}}
   \ldots \partial x_{i_{k}}}
\big|
\leq
C ( 1 + | x |^{ \chi } ),
\quad
k = 2,3,\cdots, 2p+2.
\end{equation}
By the Taylor expansion, we get 
\begin{eqnarray}
\label{eq:u-Taylor-expansion}
 && \quad
     \E \big[ u(t + h , X(t, x; t+h) ) \big]
- 
\E \big[ u ( t + h , Y(t, x; t+h) )  \big] \notag \\
 & = &
    \sum^{2p+1}_{k=1}
    \sum^{d}_{i_{1},\ldots,i_{k} = 1}
   \frac{1}{k!} 
   \frac{\partial ^{k} u ( t + h, x)}{\partial x_{i_{1}}
   \ldots \partial x_{i_{k}}}
   \E \Big[ \prod^{k}_{j=1} (\delta_{X,x})^{i_{j}}
    - \prod^{k}_{j=1} (\delta_{Y,x})^{i_{j}} \Big ]
    + \E [\mathcal{R}_{2p+2}],
\end{eqnarray}
where we used the notation \eqref{eq:notation-DeltaX-DeltaY} and denote
$$ 
\big( \delta_{X,x} \big)^{\alpha}
   := \prod_{j=1}^{2p+2}  (\delta_{X,x})^{i_j} , 
   \quad 
 \big( \delta_{Y,x} \big)^{\alpha}
 := \prod_{j=1}^{2p+2}  (\delta_{Y,x})^{i_j},
 \quad
 \alpha = ( i_1, i_2, ..., i_{2p+2} )
$$
and
\begin{eqnarray}
 \label{rem:R}
  \mathcal{R}_{2p+2} &:= & \sum_{|\alpha|=2p+2}
  \frac{|\alpha|}{\alpha!}
  \bigg( 
  \int_0^1 (1-s)^{|\alpha|-1}D^{\alpha}
   u \big( t+ h, x + s \delta_{X,x} \big) \dd s
   \, ( \delta_{X,x} )^{\alpha}  
   \notag \\
& &\qquad
     - 
  \int_0^1 (1-s)^{|\alpha|-1}D^{\alpha}
   u \big( t +h, x + s \delta_{Y,x} \big) \dd s
    \, ( \delta_{Y,x} )^{\alpha
    } 
    \bigg).
\end{eqnarray}
Combining the assumption \eqref{eq:error-of-general-one-step-approximation} with  \eqref{eq:1st-partial-u}, \eqref{eq:kth-partial-u} yields
\begin{equation} \label{eq:Taylor-expansion-estimate1}
\sum^{2p+1}_{k=1}\sum^{d}_{i_{1},\ldots,i_{k} = 1}
\frac{1}{k!} \Big | \frac{\partial ^{k} u ( t+h, x )}{\partial x_{i_{1}}
	\ldots \partial x_{i_{k}}} \Big|
	\cdot
\Big| \E \Big[ \prod^{k}_{j=1} (\delta_{X,x})^{i_{j}}
- \prod^{k}_{j=1} (\delta_{Y,x})^{i_{j}} \Big ] \Big| 
\leq
C ( 1 + | x |^{\chi} )
K ( x )
h^{p+1} .
\end{equation}
Using the H\"older inequality,
\eqref{eq:asum_solution_moment_bound}, 
\eqref{eq:asum_approximation_moment_bound}
and
\eqref{eq:kth-partial-u}
 one can derive 
\begin{eqnarray}
 \label{eq:Taylor-expansion-estimate2}
\big | \E [ \mathcal{R}_{2p+2}] \big|
  & \leq &
  \sum_{|\alpha|=2p+2} 
  C  \int_0^1 \| D^{\alpha}
   u \big( t +h,  (1-s)x + sX(t,x;t+h)\big) \|_{L^2 (\Omega; \R)} \dd s
   \| ( \delta_{X,x} )^{\alpha}
   \|_{L^2 (\Omega; \R)}  \notag \\
&& \quad +
  \sum_{|\alpha|=2p+2} 
  C  \int_0^1 \| D^{\alpha}
   u \big( t+h, (1-s)x + sY(t,x;t+h)\big) \|_{L^2 (\Omega; \R)} 
   \dd s
    \| ( \delta_{Y,x} )^{\alpha} 
   \|_{L^2 (\Omega; \R)} \notag \\
    & \leq &
    \sum_{|\alpha|=2p+2} 
    C  
    \big( 1+ | x |^{\chi} 
    + \| X(t,x;t+h) \|^{\chi}_{L^{2\chi} (\Omega; \R)} \big) 
     \,
   K ( x )
   h^{p+1} \notag 
   \\
    && \quad +
    \sum_{|\alpha|=2p+2} 
    C  
    \big( 1+ | x |^{\chi} 
    + \| Y(t,x;t+h) \|^{\chi}_{L^{2\chi} (\Omega; \R)} \big) 
    \,
    K ( x )
    h^{p+1}  \notag \\ 
& \leq &
C ( 1 + | x |^{ \beta \chi} )
K ( x )
h^{p+1},
\end{eqnarray}
where we used the assumption
	$ q_0 = 2 {\chi} $,  $q_1 = 4 p + 4$, 
	and $ q_2 = 2 {\chi} , \beta \geq 1. $
Plugging \eqref{eq:Taylor-expansion-estimate1} and 
\eqref{eq:Taylor-expansion-estimate2} into \eqref{eq:u-Taylor-expansion} enables us to get
\begin{equation}
\big | \E \big[ \varphi \big( X(t,x;t+h) \big) \big] 
-
\E \big[ \varphi \big( Y(t,x;t+h) \big) \big] \big| 
\leq 
C ( 1 + | x |^{\beta \chi} ) 
K ( x )
h^{p+1}.
\end{equation}
The proof of the lemma is thus complete.
\qed 
\section{Proofs of Lemmas in Section \ref{sec:one-step-euler}}
\label{appendix:proof_EM_one_step_error}
\textit{We first prove Lemma \ref{lem:norm-estimate-Euler-method}}. 
According to  
\eqref{eq:solution-of-SDE}
and
\eqref{eq:one_step_Euler_Maruyama_method}, 
one can use the  H\"older inequality, the moment inequality \cite{E2021applied, Zhang2017numerical}, \eqref{eq:fXt-fxs-Holder}
and \eqref{eq:gXt-gxs-Holder}, to acquire
\begin{eqnarray}
 &&
     \E [
       | \delta_{X,x}
         - \delta_{Y_E,x} |^{2q} ]   \notag \\
      & = & 
         \E \Big [ \Big |
          \int^{t+h}_t 
          f \big( X( t,x; s ) \big) 
            - f (x) 
            \, \dd s
            +
             \int^{t+h}_t g
          \big( X( t,x; s ) \big) 
            - g (x) 
            \, \dd W (s)
             \Big |^{2q} \Big ]   \notag \\
      & \leq &
           C \E \Big [ \Big|
            \int^{t+h}_t 
             f \big( X( t,x; s ) \big)
               - f (x) 
               \, \dd s 
                \Big|^{2q} \Big ]  
        +
          C   \E \Big [ \Big| 
            \int^{t+h}_t g
              \big( X( t,x; s ) \big) 
                - g (x)
                \,  \dd W ( s ) 
                  \Big|^{2q} \Big ]      \notag \\
      & \leq &
          C h^{2q-1} 
           \int^{t+h}_t 
             \big\| 
              f \big( X( t,x; s ) \big)
               - f (x)
                 \big\|_{L^{2q} (\Omega, \R^d ) }
                 ^{2q}
                 \, \dd s 
          + 
          C h^{q-1} 
           \int^{t+h}_t 
             \big \|  g \big( X( t,x; s ) \big) 
              - g(x)  \big \|_{L^{2q} (\Omega, \R^{d \times m}) }
               ^{2q}
                 \, \dd s     \notag  \\
       &  \leq &
      C (1 + |x|^{2q (4r+1) }) h^{2q},
\end{eqnarray}
which implies \eqref{eq:norm-error-one-step-Euler-method}, as required.
Similar arguments, together with \eqref{ass:drift-f}
and 
\eqref{ass:diff-g}
yield that
\begin{eqnarray}
  \label{eq:error_norm_one_step_EM}
     \E \big [ 
       | \delta_{Y_E,x} |^{2q} \big ]
  & = &
      \E \Big[ \Big| 
       \int^{t+h}_t f ( x ) 
        \, \dd s
        +
         \int^{t+h}_t g ( x ) 
          \, \dd W ( s ) 
            \Big|^{2q}  \Big ]   \notag  \\
  & \leq  &
         C h^{2q } 
           \big| 
            f ( x ) \big|^{2q} 
               +  C h^{q}  
                 \|
                  g ( x ) 
                  \|^{2q}   \notag \\
   &  \leq  &
      C (1 + |x|^{2q ( 2r+ 1)})h^{q},
\end{eqnarray}
so that \eqref{eq:norm-estimate-one-step-Euler-method} is validated.
The estimate of $ \delta_{X,x} $ as \eqref{eq:norm-estimate-one-step-exact-solution}
can be  done similarly. 
\qed
\vskip 10pt 
\textit{We then prove Lemma \ref{lem:error-one-step-Euler-method}.}
We treat the case $ s = 1 $ first.
Using the multi-dimensional It\^{o} formula,
the Cauchy-Schwarz inequality,
\eqref{ass:drift-f}
and 
\eqref{ass:diff-g}
shows
\begin{eqnarray}
\label{eq:s1-error-one-step-Euler-method}
   & &
       \Big| \E \big [ 
       ( \delta_{X,x} )^{i_1}
    -  ( \delta_{Y_E,x} )^{i_1} 
        \big ] \Big |   
    =  
      \Big| \E  \Big [ 
      \int^{t+h}_t f^{i_1}
         \big( X(s) \big)
       - f^{i_1}(x) 
        \, \dd s 
         \Big ] \Big|    \notag \\
  & = &
      \Big|  \E \Big [
     \int_t^{t+h} \int^s_t
        \left\langle
      	 D f^{i_1} \big( X(r) \big)
	    ,
	    f \big( X(r) \big) 
	   \right\rangle 
           \, \dd r \, \dd s  \Big]
    +  \E  \Big[ 
        \int_t^{t+h}\int^s_t
         \left \langle
      	  D f^{i_1} \big( X(r) \big)
	     ,
      	 g 
	      \big( X(r) \big )
	    \, \dd W (r) 
   	      \right \rangle 
         \, \dd s \Big ]   \notag \\
  & & \quad  + 
        \frac{1}{2}   
         \E \Big[ 
           \int_t^{t+h} \int^s_t
             \text{trace}
     	       \big( g^T ( X(r) )
	             D^2 f^{i_1} \big( X(r) \big)
	            g( X(r))
	           \big) \, \dd r \, \dd s \Big] \Big|   \notag \\
  &  \leq &
         \Big| \E \Big [
            \int_t^{t+h} \int^s_t
              \big| D f^{i_1} \big( X(r) \big)  \big|
	           \big| f \big( X(r) \big) \big| 
         \, \dd r \, \dd s \Big ] \Big|   \notag  \\
  &  &  \quad +  
            C  \Big| \E \Big[ 
              \int_t^{t+h} \int_t^{s}
                \big\| D^2 f^{i_1} \big( X(r) \big) \big\|
                 \big\|
                    g \big( X(r) \big)
                  \big\|^2 
                  \, \dd r  \, \dd s \Big ] \Big |  \notag \\
   &  \leq &
        C ( 1+|x|^{4r+1} ) h^{2}
     +
        C ( 1 + | x |^{ 2 r + 2 \rho - 1} ) h^{2}  \notag \\
   &  \leq &
        C ( 1+|x|^{4r+1} ) h^{2}.
\end{eqnarray}
For $ s=2,$ we first derive
\begin{eqnarray}
\label{eq:error-of-Euler-method-s-2}
    & &
      \big | \E \big [
         (\delta_{X,x})^{i_1}(\delta_{X,x})^{i_2}
      -   ( \delta_{Y_E,x})^{i_1} (\delta_{Y_E,x})^{i_2} 
            \big ] \big|   \notag  \\
    & \leq  &
       \big| \E \big [
          ((\delta_{X,x})^{i_1}
      -   (\delta_{Y_E,x})^{i_1})
           (\delta_{Y_E,x})^{i_2}  
           \big ] \big|
      +  \big| 
            \E \big [
            (\delta_{X,x})^{i_1} ((\delta_{X,x})^{i_2}
           - ( \delta_{Y_E,x})^{i_2})
              \big ] \big | 
    = :
       I_1 + I_2.
\end{eqnarray}
Before proceeding further with the estimate of $ I_{1}, $
we decompose it as follows:
\begin{eqnarray}
       I_{1}
  & = &
      \Big| \E \Big[ \Big (
        \int_t^{t+h} 
      f^{i_1} \big( X(s) \big)
      - f^{i_1}(x) 
      \, \dd s
      +  \int_t^{t+h} 
       g_{i_1 } 
       \big(
       X(s)
       \big)
      - g_{i_1} (x) 
      \, \dd W ( s )
      \Big )       \notag  \\
  & & \qquad
    \times  \Big( \int_t^{t+h} 
          f^{i_2}(x) 
       \, \dd s
     +  \int_t^{t+h}
        g_{i_2 } (x) 
        \, \dd W (s)
     \Big) \Big] \Big|    \notag \\
   & \leq & 
     \Big| \E \Big[ \Big ( \int_t^{t+h}
          f^{i_1} \big( X(s) \big)
      - f^{i_1}(x) 
      \, \dd s \Big) \cdot
         \int_t^{t+h}
         f^{i_2}(x) 
         \, \dd s 
     \Big] \Big|     \notag \\
  & & \quad
       + \Big| \E \Big[ \Big(
          \int_t^{t+h} 
          g_{i_1 } \big( X(s) \big)
       -  g_{i_1 } (x)  
       \, \dd W(s) 
       \Big) \cdot  
          \int_t^{t+h}
           f^{i_2}(x)
       \, \dd s  \Big] \Big|    \notag \\
  &  & \quad +
      \Big| \E \Big[ \Big(
          \int_t^{t+h} 
          f^{i_1} \big( X(s) \big)
       -  f^{i_1}(x) 
       \, \dd s 
      \Big) \cdot
          \int_t^{t+h} 
           g_{i_2}(x) 
           \, \dd W (s) 
      \Big] \Big|  \notag  \\
  & & \quad + 
       \Big| \E \Big[ \Big( 
          \int_t^{t+h} 
          g_{i_1} \big( X(s) \big)
        - g_{i_1} (x) 
        \, \dd W(s) 
        \Big) \cdot 
           \int_t^{t+h} 
            g_{i_2 }(x) 
            \, \dd W (s)
         \Big] \Big|   \notag \\
  & = &
      \Big| \E \Big[ \Big(
             \int_t^{t+h}
             f^{i_1}\big( X(s) \big)
        -    f^{i_1}(x) 
        \, \dd s
       \Big) \cdot 
             \int_t^{t+h}
              f^{i_2}(x) 
              \, \dd s 
       \Big] \Big|    \notag \\
  &  & \quad +
    \Big| \E \Big[ \Big( 
             \int_t^{t+h} 
              f^{i_1} \big( X(s) \big)
       -    f^{i_1}(x)
       \, \dd s  
       \Big) \cdot
             \int_t^{t+h}
              g_{i_2 }(x) 
              \, \dd W ( s ) 
        \Big] \Big|    \notag \\
  & & \quad  + 
        \Big| \E \Big[ \Big( 
            \int_t^{t+h}
            g_{i_1} \big( X(s) \big)
        -   g_{i_1 } (x) 
        \, \dd W (s) 
        \Big) \cdot 
            \int_t^{t+h}
             g_{i_2 }(x)
             \, \dd W (s)
        \Big] \Big|  \notag  \\
   & = : &
        I_{1,1} +  I_{1,2}
         +  I_{1,3}.
\end{eqnarray}
Since $ I_{1,1} $ is an easy term, 
we treat it first.
In view of \eqref{eq:f-growth} and \eqref{eq:fXt-fxs-Holder}, we get
\begin{equation}
   I_{1,1}
     \leq
    h | f^{i_2}(x) |
     \int_t^{t+h}
             \big| \E [ f^{i_1} ( X(s) )
        -    f^{i_1}(x) ]  \big|
        \, \dd s
         \leq
    C ( 1 + |x|^{6 r + 2} )  h ^{\frac{5}{2}}.
\end{equation}
Regarding $ I_{1,2}, $
with the aid of the H\"older inequality
and the It\^{o} isometry, 
one employs \eqref{ass:diff-g}
and \eqref{eq:fXt-fxs-Holder} to obtain
\begin{eqnarray}
  I_{1,2}
    &  \leq  &
           \int_t^{t+h}
           \| f^{i_1} \big( X(s) \big)
            - f^{i_1}(x)  
            \|_{L^2 ( \Omega, \R )}
             \, \dd s 
             \cdot
                \Big\| 
                \int_t^{t+h}
                 g_{i_2} (x) 
                 \, \dd W(s)
                  \Big\|_{L^2 ( \Omega, \R )}  \notag \\
    &  \leq  &
           C ( 1 + |x|^{4 r + 1} ) h^{\frac{3}{2}} 
                  \Big( 
                  \int_t^{t+h} \E \big[ 
                  | g_{ i_2 } (x) |^2 \big]
                      \,\dd s
                       \Big )^{\frac{1}{2}}   \notag \\
    & \leq &
        C (1+ |x|^{4r + \rho + 1}) h^2 .
\end{eqnarray}
Now we come to the estimate of $ I_{1,3}$, 
which requires more careful arguments.
In light of the multi-dimensional It\^{o} formula,
one can separate the considered term 
$ I_{1,3} $
as follows:
\begin{eqnarray}
\label{eq:decompose-of-I13}
    I_{1,3}
      & \leq &
     	\Big| \E \Big[ \Big( 
	     \int_t^{t+h}
	      \Big\langle 
	       \int_t^{s}
	         D g_{i_1} (X(r))
	          \cdot
	          g(X(r)) \dd W(r)
	          ,
	          \dd W(s)
	         \Big\rangle
	         \Big)
	          \int_{t}^{t+h}
	           g_{i_2} (x)
	      \, \dd W ( s ) 
	       \Big ] \Big |    \notag \\
      & & \quad +
         \Big| \E \Big[ 
           \Big( 
           \int_t^{t+h} 
             \Big \langle
             \int_t^{s}
               D g_{i_1} \big( X(r) \big)
                \cdot
                 f(X(r)) \dd r
                 ,
                 \dd W ( s )
                 \Big \rangle   \Big)
               \int_{t}^{t+h}
                 g_{i_2} (x)
             \, \dd W(s) 
               \Big ] \Big |    \notag \\
      & & \quad +  
            \frac{1}{2} 
             \sum_{l=1}^m
              \Big| \E \Big[ 
               \int_t^{t+h} \int_t^{s}
                  \text{trace} 
                    \big(
                     g^T(X(r)) D^2 g^{i_1,l}(X(r)) g(X(r))
                  \,  \dd r  
                  \, \dd W^{l}(s) 
                     \big)
\!
            \int_{t}^{t+h}
               g_{i_{2}} (x) 
              \, \dd W (s)
                \Big ] \Big |   \notag \\
      & =: &
            I^{(1)}_{1,3} + I^{(2)}_{1,3} +I^{(3)}_{1,3}. 
\end{eqnarray}
The three terms 
 $ I^{(1)}_{1,3}, $
 $ I^{(2)}_{1,3}, $
 $ I^{(3)}_{1,3}  $
will be estimated separately.
The generalized It\^{o} isometry yields that
\begin{eqnarray}
	  I^{(1)}_{1,3}  
	     &  =  &
	             \Big|  \E \Big[  \int_t^{t+h}  \Big( 
	               \int_t^{s}  D g_{i_1} \big( X(r) \big)
	                  \cdot  
	                 g \big( X(r) \big)  \dd W(r)   \Big)
	                     g_{i_2} (x)   \, \dd  s  
	                      \Big ] \Big |      \notag   \\
	       & =  &
	              \Big|   \int_{t}^{t+h}  \E    \Big[  
	                \int_{t}^{s}   D g_{i_1} \big( X(r) \big)
	                  g \big( X(r) \big)  \dd W(r)  \Big]
	                   g_{i_2} (x)    \dd  s   \Big|      \notag  \\
	     &   =   & 0.
	\end{eqnarray}
In order to properly handle
 $ I^{(2)}_{1,3}, $
we utilize the H\"older inequality,
 \eqref{ass:drift-f} and \eqref{ass:diff-g}
to obtain
\begin{eqnarray}
  \label{eq:esti-I13}
	 I^{(2)}_{1,3}
    & \leq &
        \Big \|
         \int_t^{t+h} 
          \Big\langle
            \int_t^{s}
              D g_{i_1} \big( X(r) \big)
              \cdot
               f(X(r)) \dd r
               ,
               \dd W ( s )
            \Big\rangle
            \Big \|_{L^2 ( \Omega, \R )}
          \Big \|
              \int_{t}^{t+h}
                g_{i_2} (x) 
           \, \dd W(s)
             \Big\|_{L^2 ( \Omega, \R )}      \notag \\
     &  \leq &
            C ( 1 + |x|^{2r+\rho}) h^{\frac{3}{2}}  
         \cdot
            ( 1 + |x|^{\rho}) h^{\frac{1}{2}}    \leq 
            C ( 1 + |x|^{2r+2\rho}) h^{2}.
\end{eqnarray}
For the term $ I^{(3)}_{1,3},$
similar techniques used in
\eqref{eq:esti-I13} help us to show
\begin{equation}
	 I^{(3)}_{1,3}
          \leq
               C ( 1 + |x|^{4 \rho-2} ) h^{2}.
\end{equation}
Substituting the above estimates for
$ I^{(1)}_{1,3},$
$ I^{(2)}_{1,3} $ and $ I^{(3)}_{1,3} $
into \eqref{eq:decompose-of-I13} gives
\begin{equation}
   I_{1,3}
	 \leq 
	   C ( 1 + |x|^{2r+2\rho}) h^{2}. 
\end{equation}
Collecting the estimates of the three parts
of the $ I_{1}, $ we thus arrive at
\begin{equation}
  \label{eq:esti-I1}
    I_{1}
      \leq
        C ( 1 + |x|^{ 6 r + 2} ) 
        h^{2}.
\end{equation} 
In the sequel, we focus on the estimate
of $ I_{2}.$ Likewise, we decompose it as follows:
\begin{eqnarray}
	I_{2}
	& = &
	\Big| \E \Big[ \Big (
	\int_t^{t+h} 
	f^{i_2} \big( X(s) \big)
	- f^{i_2}(x) 
	\, \dd s
	+  \int_t^{t+h} 
	g_{i_2} 
	\big(
	X(s)
	\big)
	- g_{i_2} (x) 
	\, \dd W ( s )
	\Big )       \notag  \\
	& & \qquad
	\times  \Big( \int_t^{t+h} 
	f^{i_1} \big( X(s) \big) \, \dd s
	+  \int_t^{t+h}
	g_{i_1 } \big( X(s) \big) 	\, \dd W (s)
	\Big) \Big] \Big|    \notag \\
	& \leq & 
	\Big| \E \Big[ \Big ( \int_t^{t+h}
	f^{i_2} \big( X(s) \big)
	- f^{i_2}(x) 
	\, \dd s \Big) \cdot
	\int_t^{t+h}
	f^{i_1} \big( X(s) \big) \, \dd s 
	\Big] \Big|     \notag \\
	& & \quad
	+ \Big| \E \Big[ \Big(
	\int_t^{t+h} 
	g_{i_2} \big( X(s) \big)
  -  g_{i_2} (x)  \, \dd W(s) 
	\Big) \cdot  
	\int_t^{t+h}
	f^{i_1} \big( X(s) \big) 
	\, \dd s  \Big] \Big|    \notag \\
	&  & \quad +
	\Big| \E \Big[ \Big(
	\int_t^{t+h} 
	f^{i_2} \big( X(s) \big)
	-  f^{i_2}(x)  \, \dd s 
	\Big) \cdot
	\int_t^{t+h} 
	g_{i_1} \big( X(s) \big)  \, \dd W (s) 
	\Big] \Big|  \notag  \\
	& & \quad + 
	\Big| \E \Big[ \Big( 
	\int_t^{t+h} 
	g_{i_2} \big( X(s) \big)
	- g_{i_2} (x) 	\, \dd W(s) 
	\Big) \cdot 
	\int_t^{t+h} 
	g_{i_1 } \big( X(s) \big)  \, \dd W (s)
	\Big] \Big|   \notag \\
	& = : &
	I_{2,1} +  I_{2,2} +  I_{2,3} +  I_{2,4}.
\end{eqnarray}
Next we focus on the estimations of $ I_{2,1}$-$I_{2,4}$. Similar arguments in  the estimate of $I_{1,2}$ lead to
 \begin{align}
 	 I_{2,1}   &  \leq  C ( 1 + |x|^{6 r + 2} ) h^{2}.  \\
 	 I_{2,2}   &  \leq  C ( 1 + |x|^{4 r + \rho + 1} ) h^{2}.       \\
 	 I_{2,3}   &  \leq  C ( 1 + |x|^{4 r + \rho + 1} ) h^{2}.        	 
 	\end{align}
Now we are in a position to estimate $I_{2,4}. $ Using the multi-dimensional It\^{o} formula, we can get
  \begin{eqnarray}
  	 I_{2, 4}  & \leq &
 \Big| \E \Big[ \Big( 
 \int_t^{t+h}
 \Big \langle  \int_t^{s}
 D g_{i_2} \big( X(r) \big) \cdot g \big( X(r) \big) \dd W(r), \dd W(s)
 \Big\rangle  \Big)  \!
 \int_{t}^{t+h} g_{i_1} \big( X(s) \big)
 \, \dd W ( s )  \Big ] \Big |    \notag \\
 & & \quad +
 \Big| \E \Big[ 
 \Big(  \int_t^{t+h}  \Big \langle
 \int_t^{s} D g_{i_2} \big( X(r) \big)
 \cdot f(X(r)) \dd r, \dd W ( s )
 \Big \rangle   \Big)  \!
 \int_{t}^{t+h} g_{i_1} \big( X(s) \big)
 \, \dd W(s)  \Big ] \Big |    \notag \\
 & & \quad +  
 \frac{1}{2}  \sum_{l=1}^m
 \Big| \E \Big[   \int_t^{t+h} \int_t^{s}
 \text{trace}  \big(
 g^T(X(r)) D^2 g^{i_2,l}\big( X(r) \big) g \big( X(r) \big)
 \,  \dd r  \, \dd W^{l}(s)  \big)
 \!
 \int_{t}^{t+h}
 g_{i_{1}} \big( X(s) \big)  \, \dd W (s)
 \Big ] \Big |   \notag \\
 & =: &
 I^{(1)}_{2,4} + I^{(2)}_{2,4} +I^{(3)}_{2,4}. 
  	\end{eqnarray}
Next we come to the estimate of $  I^{(1)}_{2,4},$ which requires more careful arguments. According to the generalized It\^{o} isometry, the H\"older inequality,
\eqref{ass:diff-g} and \eqref{eq:gXt-gxs-Holder} yields that
\begin{eqnarray}
  I^{(1)}_{2,4}  & = & 
	 \Big| \E \Big[ \Big(  \int_t^{t+h}
	  \Big \langle  \int_t^{s}
	  D g_{i_2} \big( X(r) \big) \cdot g \big( X(r) \big) \dd W(r),
	   \dd W(s) \Big\rangle  \Big)  \!
	  \int_{t}^{t+h} g_{i_1} \big( X(s) \big)
	   - g_{i_1} (x)  \, \dd W ( s )  \Big ] \Big |    \notag \\
	 & & \quad +
	   \Big| \E \Big[ \Big(  \int_t^{t+h}
	  \Big \langle  \int_t^{s}
	  D g_{i_2} \big( X(r) \big) \cdot g \big( X(r) \big) \dd W(r),
	  \dd W(s) \Big\rangle  \Big)  \!
	  \int_{t}^{t+h} g_{i_1} (x)  \, \dd W ( s )  \Big ] \Big |    \notag \\
     & = & 
     \Big| \int_t^{t+h}  \E \Big[ \big( 
      \int_t^{s} D g_{i_2} \big( X(r) \big) g \big( X(r) \big) 
      \dd W(r) \big)  \!
     \big( g_{i_1} \big( X(s) \big)
     - g_{i_1} (x) \big)  \, \dd s  \Big ] \Big |    \notag \\
     & \leq &
     \int_t^{t+h} 
    \big \| \int_t^{s} D g_{i_2} \big( X(r) \big)
     g \big( X(r) \big) \dd W ( r ) \big \|_{L^2 ( \Omega, \R )} 
     \cdot
     \big \| g_{i_1} \big( X(s) \big) - g_{i_1} (x)
     \big\|_{L^2 ( \Omega, \R )}    \dd s   \notag \\
     & \leq &
     C ( 1 + |x|^{ 2 r + 3 \rho -1}) h^{2}.  
	\end{eqnarray}
Similar arguments in the estimate of $ I^{(2)}_{1,3}$ lead to
\begin{align}
	  I^{(2)}_{2,4}
	   & \leq  C ( 1 + |x|^{ 2 r + 2 \rho} ) h^{2},   \\
	   I^{(3)}_{2,4}
	   & \leq  C ( 1 + |x|^{ 4 \rho -2 } ) h^{2}.   
  \end{align}
From the above estimations, it follows that
\begin{equation}
   \label{eq:esti_I_2}
	  I_2  \leq  C ( 1 + |x|^{ 6 r + 2 } ) h^{2}. 
	\end{equation}
This together with \eqref{eq:esti-I1}
helps us to derive from 
\eqref{eq:error-of-Euler-method-s-2}
that
\begin{equation}
      \big| \E \big[ 
          (\delta_{X,x})^{i_1}(\delta_{X,x})^{i_2}
           - (\delta_{Y_E,x})^{i_1} ( \delta_{Y_E,x})^{i_2}\
           \big ] \big| 
            \leq
           	C ( 1 + |x|^{6 r + 2} ) h^{2}.
\end{equation}
To handle the case $ s = 3,$
we first infer that
\begin{eqnarray}
  & &
     \big| \E \big [ 
       (\delta_{X,x})^{i_1}(\delta_{X,x})^{i_2} (\delta_{X,x})^{i_3}
          - 
        ( \delta_{Y_E,x})^{i_1} ( \delta_{Y_E,x})^{i_2} ( \delta_{Y_E,x})^{i_3} 
           \big ] \big|     \notag  \\
   & \leq &
       \big| \E \big[ 
       (\delta_{X,x})^{i_1}((\delta_{X,x})^{i_2}(\delta_{X,x})^{i_3}
      -
       ( \delta_{Y_E,x})^{i_2} ( \delta_{Y_E,x} )^{i_3})
        \big ] \big|
       + 
         \big| \E \big[
          ((\delta_{X,x})^{i_1}-(
                                \delta_{Y_E,x})^{i_1})
            ( \delta_{Y_E,x})^{i_2} ( \delta_{Y_E,x})^{i_3}      \big ] \big|  \notag \\
   & \leq &
        \big| \E \big [
          (\delta_{X,x})^{i_1} (\delta_{X,x})^{i_2}
           ((\delta_{X,x})^{i_3}- ( \delta_{Y_E,x})^{i_3})
              \big ] \big|
           + \big| \E \big [ 
           (\delta_{X,x})^{i_1} ( \delta_{Y_E,x})^{i_3}
             ( (\delta_{X,x})^{i_2}- ( \delta_{Y_E,x})^{i_2} ) \big ] \big|  \notag \\
   & & \quad +
        \big| \E \big[
         ((\delta_{X,x})^{i_1}- ( \delta_{Y_E,x})^{i_1})
            ( \delta_{Y_E,x})^{i_2} ( \delta_{Y_E,x})^{i_3} 
            \big ] \big|.
\end{eqnarray}
Combining the H\"older inequality, 
\eqref{eq:norm-error-one-step-Euler-method}, 
\eqref{eq:norm-estimate-one-step-Euler-method}, 
and 
\eqref{eq:norm-estimate-one-step-exact-solution}
shows 
\begin{equation}
    \big| \E \big[ 
      (\delta_{X,x})^{i_1} (\delta_{X,x})^{i_2} (\delta_{X,x})^{i_3}
        - 
        ( \delta_{Y_E,x})^{i_1} ( \delta_{Y_E,x})^{i_2} ( \delta_{Y_E,x})^{i_3} 
          \big ] \big|
      \leq
       C ( 1 + |x|^{8r+3} )h^{2}.
 \end{equation}
Putting these estimates with $s = 1, 2, 3$ together gives the desired assertion.
\qed
\section{Proof of Lemma \ref{lem:norm-esti-one-step-modified-EM}}
\label{appendix:moment_esti_one_step_EM_MES}
Using the H\"older inequality, \eqref{con:modified_EM_sup_lin_f}
and \eqref{con:modified_EM_sup_lin_g}
gives
\begin{eqnarray}
       \E [ | \delta_{Y_E,x}
          -
            \delta_{Y,x} |^{2q}]
     & = &
     \E \Big[ \Big|
      \int_t^{t+h} 
      f(x) - \bar{f}_{h}(x)  
        \, \dd s  
   + 
        \int_t^{t+h} 
       g(x) - \bar{g}_{h}(x)  
         \, \dd W(s) \Big|^{2q} \Big]    \notag  \\
 & \leq &
       C h^{2q-1}
        \int_t^{t+h}  
          \big\| 
          f(x) - \bar{f}_{h}(x) 
          \big\|_{L^{2q} (\Omega, \R^d ) }
                 ^{2q}
           \, \dd s    \\
  & & \quad +  
     C \big\| g(x) - \bar{g}_{h}(x) 
       \big\|_{L^{4 q} (\Omega, \R^{d \times m}) }^{2q}  \cdot
      \big\| W(t+h) - W(t) 
      \big\|_{L^{4 q} (\Omega, \R^m ) }
        ^{2q}	
     \notag \\
 & \leq &
    C (1+|x|^{2 q r_1}) h^{2q(\mathsf{a}+1)}
  + 
    C (1+|x|^{2 q \rho_1}) h^{2q(\mathsf{b}+\frac{1}{2})},
\end{eqnarray}
which implies the first assertion.
Similar arguments as above
lead to
\begin{eqnarray}
       \E \big[ \big| 
           \delta_{Y,x} \big|^{2q} 
             \big]    
   &  \leq &
      C h^{2q} 
    \big\|\bar{f}_{h}(x)\big\|_{L^{2q} (\Omega, \R^d )}^{2q}
   + 
     C \big\|\bar{g}_{h}(x) \big\|_{L^{4 q} (\Omega, \R^{d \times m}) }^{2q}  \cdot
      \big\| W(t+h) - W(t) 
      \big\|_{L^{4 q} (\Omega, \R^m ) }
        ^{2q}	
    \notag  \\
    & \leq &
         C (1 +|x|^{2q(2r+r_1+\rho_1+1)}) h^{q}.
\end{eqnarray}
This gives the second assertion and finishes the proof. \qed

\section{Proof of Lemma \ref{modified_one_error}} 
\label{appendix:moment_esti_one_step_MES_sol}
A triangle inequality yields
\begin{eqnarray}
 &&
      \bigg| \E \Big [ 
          \prod_{j=1}^s (\delta_{X,x})^{i_j}
         - 
          \prod_{j=1}^s (\delta_{Y,x})^{i_j} 
       \Big] \bigg|   \notag \\
  &  \leq  &
       \bigg| \E \Big [ 
          \prod_{j=1}^s (\delta_{X,x})^{i_j}
       - 
          \prod_{j=1}^s (\delta_{Y_E,x})^{i_j} \Big] \bigg|
       + 
       \bigg| \E \Big [ 
          \prod_{j=1}^s (\delta_{Y_E,x})^{i_j}
       - 
          \prod_{j=1}^s (\delta_{Y,x})^{i_j}
       \Big] \bigg|      \notag  \\
   &   =:  &
            A_{1}+A_{2},
            \qquad
            s=1,2,3.
\end{eqnarray}
Thanks to
\eqref{eq:error-one-step-Euler-method},
for any $s = 1,2, 3 $ we have
   \begin{equation}
   	  \label{eq:A-1-estimate}
    A_{1}  
      \leq  
        C (1+ |x|^{ 8 r + 3 })h^{2}.
    \end{equation}
For $ s=1,$ 
the second term $ A_2 $ can be estimated as follows:
\begin{eqnarray}
  \label{eq:s1-A2}
       \big| \E \big [  
           ( \delta_{Y_E,x})^{i_1}
        -  (\delta_{Y,x})^{i_1} 
        \big] \big|       
 & = &
        \bigg | \E \bigg[
           \big( 
         f^{i_1}(x) - 
         \bar{f}^{i_1}_{h}(x)
         \big) h
        +   
        \big( 
        g_{i_1}(x)
        - \bar{g}^{i_1}_{h}(x)
             \big) \big( W( t + h) - W ( t ) \big)
           \bigg ] \bigg |   \notag  \\
  & \leq &
   C (1 + |x|^{r_1 }) h^{\mathsf{a}+1}
   +
   C( 1 + |x|^{\rho_1} ) h^{\mathsf{b}+1}
   \notag  \\
 & \leq &
    C (1 + |x|^{r_1 + \rho_1}) h^{ (\mathsf{a}+1)\land (\mathsf{b}+1) } .
\end{eqnarray}
For $ s=2, 3 $, with the H\"older inequality, \eqref{eq:norm-estimate-one-step-Euler-method}
and Lemma \ref{lem:norm-esti-one-step-modified-EM} at hand, one can readily get
\begin{align}
  \label{eq:s2-A2}
  \bigg| \E \bigg[ 
         \prod_{j=1}^2 (\delta_{Y_E,x})^{i_j}
     - 
          \prod_{j=1}^2 (\delta_{Y,x})^{i_j} 
     \bigg] \bigg| 
   &  \leq
        C (1+ |x|^{2r+2\rho_1+2r_1+1}) h^{(\mathsf{a}+\frac{3}{2}) \land (\mathsf{b}+1)},      \\
    \label{eq:s3-A2}
    \bigg| \E \bigg[ 
         \prod_{j=1}^3 (\delta_{Y_E,x})^{i_j}
     - 
          \prod_{j=1}^3 (\delta_{Y,x})^{i_j} 
     \bigg] \bigg|  
     &  \leq
        C (1+ |x|^{4r+3\rho_1+3r_1+2}) h^{(\mathsf{a}+2) \land (\mathsf{b}+\frac{3}{2})}. 
\end{align}
This together with
\eqref{eq:s1-A2},
\eqref{eq:s2-A2} gives
\begin{equation}
       A_2 
         \leq
            C (1+ |x|^{4r+3\rho_1+3r_1+2}) 
            h^{(\mathsf{a}+1)\land (\mathsf{b}+1)} ,
             \,\,
             s = 1,2,3.
\end{equation}
This combined with \eqref{eq:A-1-estimate} promises
\eqref{sheme:modified_EM_general_one_step_error}.
Using the H\"{o}lder inequality and 
\eqref{eq:norm-estimate-one-step-exact-solution} yields
\begin{equation}
     	\bigg\| \prod_{j=1}^{4}  
     	   ( \delta_{X, x} )^{i_j} 
     	\bigg\|_{L^2 (\Omega; \R)}
      \leq   
       	 \prod_{j=1}^{4}
       	   \Big\|   
       	   ( \delta_{X, x} )^{i_j} 
       	   \Big\|_{L^8 (\Omega; \R)}
      \leq
         C (1+|x|^{8r+4}) h^2 .
\end{equation}
Armed with  \eqref{eq:norm-esti-one-step-TEM}, 
we can derive \eqref{eq:s4_norm_esti_one_step_modified_EM} in a similar way.
\qed
\end{document}